\documentclass[10pt]{article}
\usepackage[margin=1.35in]{geometry}

\usepackage{amsmath}
\usepackage{graphicx}
\usepackage{psfrag}
\usepackage{amssymb}
\usepackage{verbatim,booktabs}
\usepackage{url}
\usepackage{mathtools}
\usepackage{float}
\usepackage{mathrsfs}
\usepackage{enumitem} 

\usepackage[usenames,dvipsnames]{pstricks}
\usepackage{epsfig}
\usepackage{pst-grad} 
\usepackage{pst-plot} 
\usepackage{xspace}
\usepackage{caption}
\usepackage{subcaption}

\usepackage{comment}
\usepackage{tikz}
\usepackage{accents}
\usepackage{multirow}
\usepackage{array}
\usepackage{longtable}
\usepackage{pdflscape}
\usepackage[numbers]{natbib}
\usepackage{hyperref}

\def\ie{{i.e.,} }
\def\eg{{e.g., }}

\newcommand{\orgdiv}[1]{#1}%
\newcommand{\orgname}[1]{#1}%
\newcommand{\orgaddress}[1]{#1}%
\newcommand{\postcode}[1]{#1}%
\newcommand{\city}[1]{#1}%
\newcommand{\country}[1]{#1}%
\usepackage[createShortEnv]{proof-at-the-end}
\newcommand{\fnm}[1]{#1}%
\newcommand{\sur}[1]{#1}%
\newcommand{\sortas}[1]{}

\newcolumntype{L}{>{$}l<{$}}
\newcolumntype{C}{>{$}c<{$}}
\newcolumntype{R}{>{$}r<{$}}

\newcommand{\bR}{\mathbb{R}}

\newcommand{\bZ}{\mathbb{Z}}

\newcommand{\cG}{{\mathcal G}}
\newcommand{\cC}{{\mathcal C}}
\newcommand{\cY}{{\mathcal Y}}

\newcommand{\cD}{{\mathcal{D}}}

\newcommand{\cP}{{\mathcal P}}
\newcommand{\cR}{{\mathcal R}}

\newcommand{\cQ}{{\mathcal Q}}
\newcommand{\cF}{{\mathcal F}}
\newcommand{\cS}{{\mathcal S}}

\newcommand{\cU}{{\mathcal U}}
\newcommand{\cV}{{\mathcal V}}
\newcommand{\cX}{{\mathcal X}}

\newcommand{\cSt}{\cS_{\mathrm{st}}}
\newcommand{\bcSt}{\bar{\cS}_{\mathrm{st}}}
\newcommand{\bcStt}[1]{\bar{\cS}^{#1}_{\mathrm{st}}}
\newcommand{\cSl}{\cS_{\mathrm{lift}}}
\newcommand{\deq}{\coloneqq}

\newcommand{\lin}[2]{{\Xi^{#1}_{#2}}}

\newcommand{\relx}[1]{\tilde{#1}}
\newcommand{\norm}[1]{{\lVert#1\rVert}}

\DeclareMathOperator{\suc}{s.t.}
\DeclareMathOperator{\proj}{proj}

\DeclareMathOperator{\dom}{dom}

\DeclareMathOperator{\rang}{range}
\DeclareMathOperator{\inter}{int}

\DeclareMathOperator{\conv}{conv}

\DeclareMathOperator{\conve}{convenv}
\DeclareMathOperator{\epi}{epi}
\DeclareMathOperator{\hyp}{hypo}
\DeclareMathOperator{\gr}{graph}
\DeclareMathOperator{\cl}{cl}
\DeclareMathOperator{\bd}{bd}

\newcommand{\scip}{$\texttt{SCIP}$\xspace}
\newcommand{\couenne}{$\texttt{Couenne}$\xspace}
\newcommand{\baron}{$\texttt{BARON}$\xspace}
\newcommand{\antigone}{$\texttt{ANTIGONE}$\xspace}
\newcommand{\miso}{$\texttt{MISO}$\xspace}
\newcommand{\minlplib}{MINLPLib\xspace}
\newcommand{\disable}{\texttt{disable}\xspace}
\newcommand{\ic}{\texttt{ic}\xspace}
\newcommand{\oc}{\texttt{oc}\xspace}
\newcommand{\oic}{\texttt{oic}\xspace}

\newcommand{\notemi}[1]{}

\usepackage{comment}

\usepackage{amsthm}
\newtheorem{theorem}{Theorem}
\newtheorem{proposition}{Proposition}
\newtheorem{lemma}{Lemma}
\newtheorem{corollary}{Corollary}

\newtheorem{definition}{Definition}

\usepackage[capitalize,noabbrev]{cleveref}

\newtheorem{example}{Example}
\newtheorem{observation}{Observation}
\Crefname{chapter}{Chap.}{Chaps.}
\Crefname{section}{Sec.}{Secs.}
\Crefname{subsection}{Subsec.}{Subsecs.}
\Crefname{proposition}{Prop.}{Props.}
\Crefname{theorem}{Thm.}{Thms.}
\Crefname{definition}{Defn.}{Defns.}
\Crefname{corollary}{Cor.}{Cors.}
\Crefname{figure}{Fig.}{Figs.}
\Crefname{observation}{Obs}{Obss.}
\Crefname{chapter}{Chap.}{Chaps.}

\newcommand{\kywrds}{global optimization, signomial programming, extended formulation, cutting plane, intersection cut, convex relaxation}

\begin{document}

\title{Cutting planes for signomial programming}

\author{\fnm{Liding} \sur{Xu}\thanks{\orgdiv{LIX CNRS}, \orgname{\'Ecole Polytechnique, Institut Polytechnique de Paris}, \orgaddress{\city{Palaiseau}, \postcode{91128}, \country{France}}.
         	E-mail: \texttt{lidingxu.ac@gmail.com}, \texttt{dambrosio@lix.polytechnique.fr}, \texttt{liberti@lix.polytechnique.fr}, \texttt{sonia.vanier@lix.polytechnique.fr}}
         	\and
\fnm{Claudia} \sur{D'Ambrosio}\footnotemark[1]
         	\and
\fnm{Leo} \sur{Liberti}\footnotemark[1]
         	\and
\fnm{Sonia} \sur{Haddad-Vanier}\footnotemark[1] \thanks{ \orgname{Universit\'e Paris 1 Panthéon-Sorbonne}, \orgaddress{\city{Paris}, \postcode{75005}, \country{France}.}}
}

\date{\today}

\maketitle

\begin{abstract}
Cutting planes are of crucial importance when solving nonconvex nonlinear programs to global optimality, for example using the spatial branch-and-bound algorithms. In this paper, we discuss the generation of cutting planes for signomial programming. Many global optimization algorithms lift signomial programs into an extended formulation such that these algorithms can construct relaxations of the signomial program by outer approximations of the lifted set encoding nonconvex signomial-term sets, i.e., hypographs, or epigraphs of signomial terms. We show that any signomial-term set can be transformed into the subset of the difference of two concave power functions, from which we derive two kinds of valid linear inequalities. Intersection cuts are constructed using signomial term-free sets which do not contain any point of the signomial-term set in their interior. We show that these signomial term-free sets are maximal in the nonnegative orthant, and use them to derive intersection sets. We then convexify a concave power function in the reformulation of the signomial-term set, resulting in a convex set containing the signomial-term set. This convex outer approximation is constructed in an extended space, and we separate a class of valid linear inequalities by projection from this approximation. We implement the valid inequalities in a global optimization solver and test them on MINLPLib instances. Our results show that both types of valid inequalities provide comparable reductions in running time, number of search nodes, and duality gap.
\end{abstract}



\emph{Key words:} \kywrds

\section{Introduction}

General nonconvex nonlinear programming (NLP) problems typically admit the following  formulation:
 \begin{equation}
\label{minlp}
	\min_{x \in \bR^n} \;  c\cdot x  \quad \suc \quad  Ax +  B g(x)\le d,
\end{equation}
where $c \in \bR^n, A \in \bR^{m \times n},  B \in \bR^{m\times \ell}, g: \bR^n \to \bR^\ell,  d \in \bR^m$. 

The mapping $g(x)$ represents a vector $(g_1(x),\dots,g_\ell(x))$ of nonconvex functions on $x$, and we denote $g_i$ as their \emph{terms}. Note that the objective function is supposed to be linear, w.l.o.g., since we can always reformulate a problem with a nonlinear objective function as the problem \eqref{minlp} above (epigraphic reformulation).

 General-purpose global optimization solvers, such as $\baron$ \cite{tawarmalani2005polyhedral}, $\couenne$ \cite{couenne}, and $\scip$ \cite{bestuzheva2023global}, are capable of solving the problem \eqref{minlp} within an $\epsilon$-global optimality. They achieve this by employing the spatial branch-and-bound (sBB) algorithm, which explores the feasible region of \eqref{minlp} implicitly, but systematically.
The sBB algorithm effectively prunes out unpromising search regions by comparing the cost of the best feasible solution found with the cost bounds associated with those regions. These cost bounds can be computed by solving convex relaxations of  \eqref{minlp}.

The backend convex relaxation algorithms implemented in many general-purpose  solvers, including \baron, \couenne,  and \scip, are linear programming relaxations. These solvers take advantage of the separability introduced in the rows of $Ax+Bg(x)$, allowing them to relax and linearize  nonlinear terms $g_i$ individually.
In the solvers' data structures, the problem \eqref{minlp} is transformed into an extended formulation:
 \begin{equation}
 \label{minlp2}
	\min_{ (x,y) \in \bR^{n+\ell}} \;  c\cdot x  \quad \suc \quad  Ax +  B y\le d \; \land \; y = g(x).
\end{equation}
All the nonlinear terms are grouped within the nonconvex constraints $y = g(x)$. These constraints give rise to a nonconvex \emph{lifted set} defined as:
\begin{equation}
\label{eq.lift}
 \cSl \deq \{(x,y) \in \bR^{n+\ell}: y = g(x)\}.	 
\end{equation}

The relaxation algorithms used by these solvers are based on factorable programming \cite{leoinbook, mccormick1976computability}: This approach treats the multivariate nonlinear terms $g_i$ as composite functions. These  algorithms typically factorize each $g_i$ into sums and products of a collection of univariate functions. If convex and concave relaxations of those univariate functions are available, these algorithms can linearize these relaxations, and yield a linear relaxation for Eq.~\eqref{minlp}. Common lists of such univariate functions, that are usually available to all sBB solvers, include $t^a$ (for $a\in\mathbb{N}$), $\frac{1}{t}$, $\log{t}$, $\exp{t}$. Some solvers also offer a choice of trigonometric functions, e.g.~\couenne. 

Most sBB solvers can handle \emph{signomial term} $\psi_{\alpha}(x) \deq x^{\alpha} =\prod_{j \in [n]} x_j^{\alpha_j}$, where  the exponent vector $\alpha$ is in $\bR^n$, but in a way that yields poor relaxations (more about this below). 
In this paper, we provide a deeper treatment of the signomial term  w.r.t. convexification and linearization within an sBB algorithm.

When all the terms in $g$ are signomial terms, the problem \eqref{minlp} falls under the category of signomial programming (SP). In this scenario, we refer to \eqref{minlp} as the \emph{natural formulation} of SP. The left-hand sides of the constraints in this formulation are referred to as \emph{signomial functions}. The lifted set $\cSl$ in the extended formulation \eqref{minlp2} is called a \emph{signomial lift}. 

Since negative entries may present in the exponent vector $\alpha$, in general, variables of SP are assumed to be positive.  The point of restriction on  SP over positive variables is simply to make the theoretical treatment more readable and
streamlined. We remark that the  techniques in this paper can also treat signomial terms in general mixed-integer NLP problems. 

In the case of SP, LP relaxations can be derived from polyhedral outer approximations of the signomial lift in its extended formulation.
A typical relaxation algorithm for SP involves factorizing the signomial term $\psi_{\alpha}(x)$ into the product of $n$ univariate signomial terms $x_i^{\alpha_i}$. After the factorization, the algorithm proceeds to convexify and linearize the intermediate  multilinear term and univariate functions. However, this factorable programming approach can lead to weak LP relaxation and introduce additional auxiliary variables that represent intermediate functions. These problems have already been discussed in the context of pure multilinear terms \cite{cafieri2010convex, costa2012relaxations, speakman2017quantifying}.

We propose two cutting plane-based relaxation algorithms for SP. In contrast to the conventional factorable programming approach, our method uses a novel reformulation of the signomial lift. We transform each nonlinear equality constraint $y_i = g_i(x)$ in \eqref{eq.lift} to an equivalent constraint $ \psi_{\beta}(u) - \psi_{\gamma}(v) = 0$, where $\beta > 0, \gamma>0$, $\max(\norm{\beta}_1, \norm{\gamma}_1)=1$, $u,v$ are  sub-vectors partitioned from $(x,y)$, and $ \psi_{\beta} , \psi_{\gamma} $ are concave functions. Thus, the nonlinear equality constraint is equivalent to two inequality constraints: $\psi_{\beta}(u) - \psi_{\gamma}(v) \le 0$ and $\psi_{\beta}(u) - \psi_{\gamma}(v) \ge 0$, with $u \in \bR_{+}^h, v \in \bR_{+}^k$ being reassignments of $(x,y)$. Our algorithms aim at generating convex relaxations of these two inequality constraints. Due to the symmetry of these two constraints, we consider convex relaxations for the first one. This reduction motivates us to construct linear valid inequalities for the nonconvex \emph{signomial-term set}:
\begin{equation}
\label{eq.sform}
	\cSt\deq\{(u, v) \in \bR_{+}^{h + k}:\, \psi_{\beta}(u) - \psi_{\gamma}(v) \le 0\},
\end{equation}
where the subscript $\mathrm{st}$ is an abbreviation for ``signomial term''.

Our first cutting plane  algorithm is based on the intersection cut paradigm \cite{conforti2011}. As shown in \Cref{sec.prem}, one can approximate a nonconvex set $\cS$ using its polyhedral outer approximation. This requires the construction of $\cS$-free sets, i.e., closed convex sets containing none of the interiors of $\cS$. The main insight about $\cS$-free sets for a nonconvex set $\cS$ is that they provide an explicit and useful description of the convex parts of the complement of $\cS$. In \Cref{sec.max} we extend several general results from the literature on maximal $\cS$-free sets. In \Cref{sec.icforsp} we give the transformation procedure leading to $\cSt$ and construct $\cSt$-free sets from the transformation. We show that these sets are also signomial-lift-free and maximal in the nonnegative orthant. We also discuss the separation of intersection cuts. 

To ensure convergence of the sBB algorithm, a common assumption for SP is that all variables are bounded. Our second cutting plane algorithm aims to approximate $\cSt$ within a hypercube. In \Cref{sec.outerforsp}, we provide an extended formulation for the convex envelope of the concave function $\psi_{\beta}$ over the hypercube. This formulation yields a convex set including $\cSt$ (which is a convex outer approximation of $\cSt$), so that we can generate outer approximation cuts by projection. We prove that $\psi_{\beta}$ is a supermodular function. For $h=2$ we provide a closed expression for its convex envelope by exploiting supermodularity, which allows us to get rid of the projection step.

For the computational part of this study, we note that signomials are one of the four main types of nonlinearities found in the mixed-integer NLP library (\minlplib) \cite{minlpsurveyAN, minlplib}. Our relaxation approach does not require factorization or the introduction of intermediate functions, so implementing the proposed cutting planes in the general-purpose solver \scip is straightforward, and the outer approximation cut algorithm is integrated in \scip since version 9.0 \cite{BolusaniEtal2024}. In \Cref{sec.comp},
we perform computational tests with instances from \minlplib and observe improvements to \scip default settings due to the proposed valid inequalities.

\subsection{Related works}

The majority of relaxations for SP are derived from its generalized geometric programming (GGP) formulation,  which is an exponential transformation \cite{duffin1970linearizing}  of its natural formulation. The exponential transformation replaces positive variables $x$ by exponentials $\exp(z)$, where $z$ are real variables. The authors of  \cite{maranas1997global} show that signomial functions in GGP are difference-of-convex (DC) functions. For the signomial function in each constraint of GGP, they  construct linear underestimators of its concave part; the author of \cite{shen2005linearization} constructs linear underestimators of the whole function via the mean value theorem. The author of \cite{xu2014global} proposes inner approximations  of GGP via the inequality of  arithmetic and geometric means (AM-GM inequality). The authors of \cite{chandrasekaran2016,dressler2022algebraic,murray2021signomial} construct  non-negativity certificates for signomial functions via the AM-GM inequality, and propose a hierarchy of convex relaxations for GGP. Exponential transformations can be combined with other variable transformations, such as  power transformations, and the inverse transformations can be approximated by piece-wise linear functions, see \cite{lin2012range,lundell2013reformulation,lundell2018}.

The solvers \scip \cite{bestuzheva2023global}, \baron \cite{tawarmalani2005polyhedral}, \antigone \cite{Misener2014}, and \miso \cite{Misener2014miso} are able to solve the natural formulation of SP or its extended formulation within a global $\epsilon$-optimality using the sBB algorithm. More precisely, \miso is a specialized solver for SP, which uses exponential transformations of some signomial terms only when necessary. For the following reasons, exponential transformations can complicate general-purpose solvers. First, in certain NLP problems, signomial terms may appear only as a subset of the nonlinear terms of $g(x)$. In such cases, solvers may need to force the inverse transformation $x_j = \ln(z_j)$, which requires additional processing for convexification algorithms. Second, when dealing with mixed-integer SP and some variables of $x$ are integer, exponential transformations cause certain components of $z$ to become discrete but not necessarily integer. As a result, the sBB algorithm must adjust its branching rules.

While much attention has been paid to the construction of relaxations for GGP, the literature on relaxations for the extended natural formulation of SP is relatively limited. The convex relaxations used in the aforementioned solvers rely mainly on factorable programming \cite{leonelson, mccormick1976computability}. Since exponential transformations are nonlinear variable transformations, it is impossible to apply the relaxations developed for the GGP formulation directly to the natural formulation.

Numerous research efforts have been devoted to improving relaxation techniques for multilinear terms and univariate/bivariate functions commonly used in factorable programming \cite{bao2015global}. Multilinear terms over the unit hypercube are vertex polyhedral and their envelopes over the unit hypercube admits simple extended formulations \cite{rikun1997convex}. In particular, there are closed forms for the convex envelopes of bilinear functions \cite{al1983jointly, mccormick1976computability} and trilinear functions \cite{meyer2004trilinear, meyer20042} over hypercubes. In \cite{sherali1997convex}, the author presents convex envelopes for multilinear functions (sum of multilinear terms) over the unit hypercube and specific discrete sets. For a comprehensive analysis of multilinear term factorization via bilinear terms, we refer to \cite{luedtke2012some, speakman2017quantifying}. Additionally, \cite{cafieri2010convex} offers an in-depth examination of quadrilinear function factorization through bilinear and trilinear terms, while \cite{costa2012relaxations} presents a computational study on extended formulations.

Convexifying univariate/bivariate functions plays an important role in the field of global optimization. In \cite{liberti2003convex}, convex envelopes for monomials with odd degrees are derived. An approach presented in \cite{locatelli2014convex} enables the evaluation of the convex envelope of a bivariate function over a polytope and separating its supporting hyperplane by solving low-dimensional convex optimization problems. The convex optimization problems are further reduced by solving a Karush-Kuhn-Tucker system \cite{locatelli2018convex}.  In \cite{locatelli2016polyhedral}, convex envelopes for bilinear, fractional, and other bivariate functions over a polytope are constructed using a polyhedral subdivision technique. The relation between
triangulation and envelope construction has been observed in \cite{tawarmalani13}, and  we refer to \cite{bao2015global,bao09}  computational studies on  triangulation-based convexification of nonconvex quadratic and multilinear terms. Additionally, \cite{nguyen2018deriving} employ polyhedral subdivision and lift-project methods to derive explicit forms of convex envelopes for various nonconvex functions, including a specific subclass of bivariate signomial terms. We refer to \cite{bure2017,kilincc2023conic,ki2015} for results on convexification of sets involving mixed-integer convex cones, as these works on convexification of such sets share some common techniques with convexification of nonconvex functions.

Convexifying high-order multivariate functions is a major challenge, and the available literature on convex underestimators for trivariate functions is relatively few. For supermodular functions, there are several classes of valid inequalities for their convex envelopes, see \cite{ahmed2011maximizing,atamturk2022supermodularity2,han2022fractional,Nemhauser1988}. In \cite{he2021new, he2022tractable}, the authors propose a novel framework for relaxing composite functions in nonlinear programs. Another approach is to use the intersection cut paradigm \cite{conforti2011} to approximate nonconvex functions. This paradigm can generate cutting planes to strengthen LP relaxations of NLP problems. Constructing intersection cuts involves finding an $\cS$-free set, where $\cS$ represents a nonconvex set defined by nonconvex functions.
The study of intersection cuts originated in the context of NLP \cite{thy1985}. Gomory later introduced the concept of corner polyhedron \cite{gomory1969some}, and intersection cuts were explored in the field of integer programming \cite{balas1971}. The modern definition of intersection cuts for arbitrary sets $\cS$ is from \cite{dey2008, glover1973}. For more comprehensive details, we refer to \cite{andersen2010, basu2019, cornuejols2015sufficiency, del2012relaxations, dey2008, richard2010group}. Recent research has revealed $\cS$-free sets for various nonconvex sets encountered in structured NLP problems. Examples include outer product sets \cite{bienstock2019}, sublevel sets of DC functions \cite{serrano2019}, quadratic sets \cite{munoz2020}, and graphs of bilinear terms \cite{fischetti2020}. Intersection cuts have also been developed for convex mixed-integer NLP problems \cite{andersen2007, belotti2015conic, kilinc-karzan2016, modaresi2016} and for bilevel programming \cite{fischetti2018}. 

\subsection{Notation}
We follow standard notation in most cases. Let $[n_1:n_2]$ stand for $\{n_1,\hdots,n_2\}$, and let $[n]$ stand for $[1:n]$.  For a vector $x \in \bR^n$, $x_j$ denotes the $j$-th entry of $x$; given $J \subseteq [n]$, $x_J = (x_j)_{j \in J}$ denotes  the sub-vector  formed by entries indexed by $J$. $\lVert \cdot \rVert_p$ denotes the $L_p$-norm ($1\le p\le +\infty$).  For  a set $X \subseteq \bR^n$, $\conv(X)$,   $\cl(X)$, $\inter(X)$, $\bd(X)$, $|X|$, $X^c$ denote the convex hull,  closure, interior, boundary, cardinality, and complement of $X$, respectively.   For a function $f$,  $\dom(f)$ and $\rang(f)$  denote the domain and range of $f$, respectively;  $\gr(f)$ denotes its graph $\{(x,t) \in \bR^{n+1}: f(x) = t\}$, $\epi(f)$ denotes its epigraph $\{(x,t) \in \bR^{n+1}: f(x) \le t\}$, and $\hyp(f)$ denotes its hypograph $\{(x,t) \in \bR^{n+1}: f(x) \ge t\}$;  if $f$ is differentiable, for a $\relx{x} \in \dom(f)$, $\nabla{f}(\relx{x})$  denotes the gradient of $f$ at $\relx{x}$ and  
\begin{equation}
    \lin{f}{\relx{x}}(x) := f(\relx{x}) + \nabla f(\relx{x}) \cdot(x - \relx{x}). \label{def.Xi}
\end{equation}
The word \textit{linearization} involves the replacement of a nonlinear function by its affine underestimators or overestimators. For example, the affine underestimators of convex functions $f$ are given as $\lin{f}{\relx{x}}(x)$ for some $\relx{x}$.

\section{Preliminaries}
\label{sec.prem}
In this section we present an overview of $\cS$-free sets and intersection cut theory. The process of constructing intersection cuts involves two fundamental steps \cite{conforti2014}: constructing $\cS$-free sets and deriving cutting planes from these sets. Since maximal $\cS$-free sets yield tightest cutting planes, one can include an optional step to check the maximality of $\cS$-free sets.

\begin{definition}
\label{def.free}
Given  a set $\cS \subsetneq \bR^p$, a closed set $\cC$ is (convex) $\cS$-free if $\cC$ is convex and $\inter(\cC) \cap \cS = \varnothing$.
\end{definition}

To construct an intersection cut, an essential requirement is the availability of a translated simplicial cone $\cR$ that satisfies two conditions: (i) $\cR$ is generated by linearly independent vectors, (ii) $\cR$ contains $\cS$, and (iii) the vertex $\relx{z}$ of $\cR$ does not belong to $\cS$.

\begin{figure}
 \centering
  \subfloat[$\cC$ as an inner approximation of $\cl(\cS^c)$.]{\includegraphics[width=0.3\textwidth]{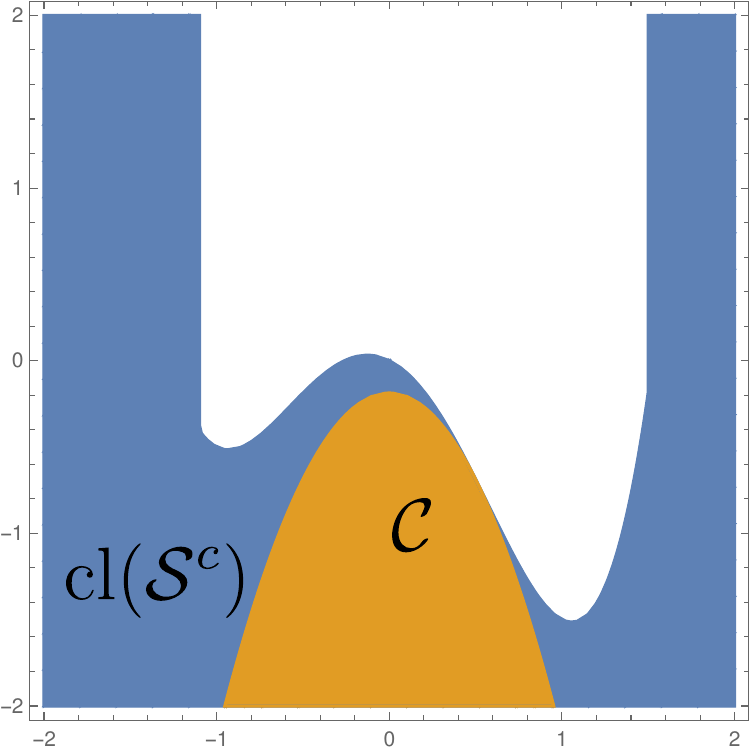} \label{fig.sdef1}}
  \hfill
  \subfloat[$\cC$ as an $\cS$-free set.]{\includegraphics[width=0.3\textwidth]{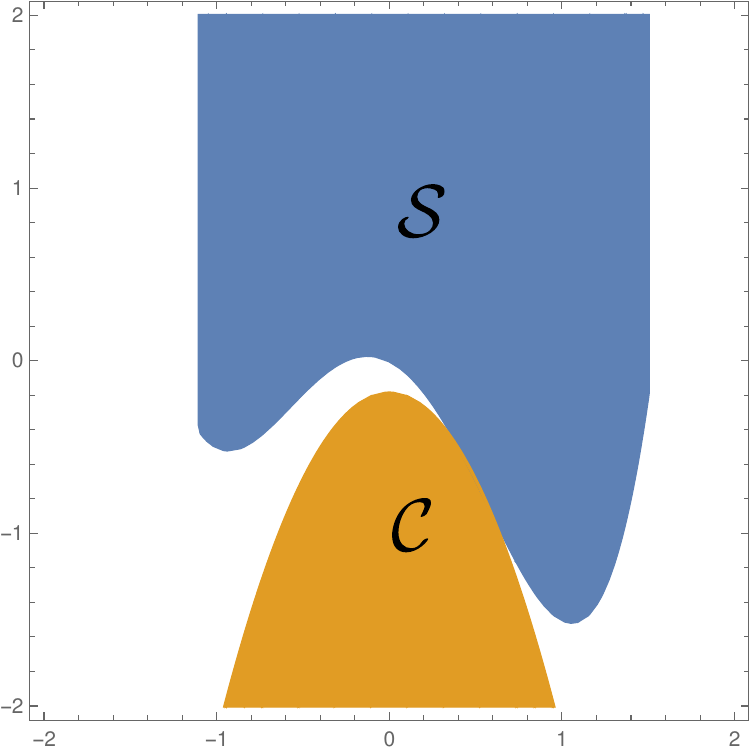}  \label{fig.sdef2}}
  \hfill
\subfloat[Simplicial cone $\cR$ and the intersection cut.]{\includegraphics[width=0.3\textwidth]{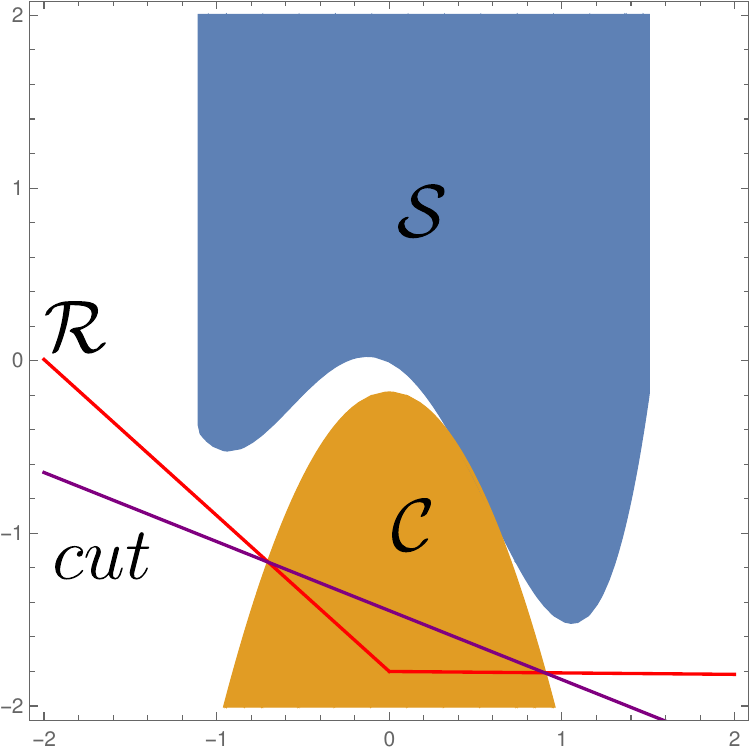}  \label{fig.sdef3}}
  \caption{An $\cS$-free set $\cC$, simplicial cone $\cR$, and intersection cut.}
  \label{fig.sdef}
\end{figure}

\Cref{fig.sdef1,fig.sdef2,fig.sdef3} give an example procedure to construct an $\cS$-free set $\cC$ and an intersection cut:  in \Cref{fig.sdef1}; we find a convex inner approximation $\cC$ of  $\cl(\cS^c)$; and we visualize the $\cS$-freeness of $\cC$ in \Cref{fig.sdef2}; then,  in \Cref{fig.sdef3}, a simplicial conic outer approximation $\cR$ of $\cS$ is used to define the intersection cut.

 We assume that  $\cR$ admits a hyper-plane representation $
 	\{z \in \bR^p: B(z - \relx{z}) \le 0\}, $
where $B \in \bR^{p \times p}$ is an  invertible matrix. For every $j\in[p]$, let $r^j$ denote the $j$-th column of $-B^{-1}$, then $r^j$ turns out to be an extreme ray of $\cR$. Thereby, $\cR$ also admits a ray representation $
   \{z \in \bR^p:\exists \mu \in \bR^p_+ \textup{ }  z =\relx{z} +  \sum_{j =1}^{p} \mu_jr^j \}. $
For every $j \in [p]$, we define the \textit{step length} from $\relx{z}$ along ray $r_j$ to the boundary  $\bd(\cC)$ as
\begin{equation}
 \label{eq.iccoef}
   \mu_j^\ast \deq \sup_{\mu_j \in [0,+\infty]}\{\mu_j: \relx{z} + \mu_j r^j \in \cC\}.
\end{equation}
Then, an intersection cut admits the form
\begin{equation}
 \label{eq.ic}
 	\sum_{j =1}^{p}  B_j (z - \relx{z})/ \mu_j^\ast\le -1,
 \end{equation}
where $B_j$ is the $j$-th row of $B$. When all step lengths are positive, the above linear inequality cuts off $\relx{z}$ from $\cS$, see for an example of an intersection cut in \Cref{fig.sdef3}.

We can obtain the sets $\cC$, $\cR$ and the vertex $\relx{z}$ by the following procedure. Suppose that we have an LP relaxation $\min_{z \in \cP} c \cdot z$ of an SP problem, where $\cP$ is a polyhedral outer approximation of the feasible set of the SP problem. If the solution to the LP problem turns out to be infeasible for the SP problem, it means that the solution does not belong to the signomial lift. In such cases, we can set $\relx{z}$ as the solution obtained from LP and let $\cC$ be the signomial-lift-free ($\cSl$-free) set. Moreover, we can extract the cone $\cR$ from the optimal LP basis defining $\relx{z}$, see \cite{conforti2014}.

One focus of our study is the construction of (maximal) $\cS$-free sets. The importance of finding \emph{maximal} sets follows from the fact that if we have two $\cS$-free sets called $\cC$ and $\cC^\ast$, where $\cC$ is a subset of $\cC^\ast$, then the intersection cut derived from $\cC^\ast$ dominates the cut derived from $\cC$ (see \cite[Remark 3.2]{conforti2011}). To give a precise characterization, we present a formal definition of maximal $\cS$-free sets.
\begin{definition}
\label{def.max}
Given  a closed convex set $\cG \subseteq \bR^p$ such that $\cS \subsetneq \cG$, an  $\cS$-free set $\cC$ is (inclusion-wise) maximal in $\cG$,  if there is no other $\cS$-free set $\cC'$ such that $\cC \cap \cG \subsetneq \cC' \cap \cG $.
\end{definition}

The above definition provides a generalization of the conventional concept of maximal $\cS$-free sets, which is a special case when $\cG = \bR^p$. Studying maximality for $\cS$-free sets in $\bR^p$ can be challenging in certain scenarios. However, \Cref{def.max} allows us to examine the intersections of $\cS$-free sets within the ground set $\cG$. This constraint is essential for our analysis, especially considering that all variables in SP are non-negative.

Next, we show how to construct $\cS$-free sets from  ``reverse'' representations of sets defined by a particular type of nonconvex functions. A function $f$ is said to be  difference-of-concave (DCC) if  there exist two concave functions $f_1,f_2$ such that $f=f_1 - f_2$.
Any DCC function is also a difference-of-convex (DC) function, and vice versa. We call a nonconvex set a \emph{DCC set}, if it admits a \emph{DCC formulation}, meaning that it is defined by a non-negative/non-positive constraint on a DCC function.  By using the \emph{reverse-minorization} technique, the following lemma provides a collection of $\cS$-free sets for DCC sets.

\begin{lemma}\cite[Prop. 6]{serrano2019}
\label{cor.dc}
Let $\cS \deq \{z \in \bR^p: f_1(z) - f_2(z) \le 0\}$, where  $f_1, f_2$ are concave functions over $\bR^p$. Then, for any $\relx{z} \in \bR^p$, $\cC \deq \{z \in \bR^p: f_1(z) - \lin{f_2}{\relx{z}}(z) \ge 0\}$ is $\cS$-free. Moreover, if $\relx{z} \in \bR^p \smallsetminus \cS$, $\relx{z} \in \inter(\cC)$.
\end{lemma}

The reverse-minorization technique involves reversing the inequality that defines $\cS$ and linearizing its convex component $-f_2$ to $-\lin{f_2}{\relx{z}}(z)$. Thus, the function $f_1(z) - \lin{f_2}{\relx{z}}(z)$ minorizes $ f_1(z) - f_2(z)$ at any $z$. The point $\relx{z}$ is referred to as the \emph{linearization point}. It is important to note that, when the shared domain $\cG$ of $f_1$ and $f_2$ is not the entire space $\mathbb{R}^p$, the set $\cS$ needs to be constrained to the \emph{ground set} $\cG$. This restriction ensures the applicability of the lemma.

\section{General results on maximality}
\label{sec.max}
In this section, we present two results on the maximality of $\cS$-free sets arising in general nonconvex NLP problems. The results are used to construct maximal signomial-lift-free sets in non-negative orthants.

\subsection{Lifted sets}
\label{sec.lift}
We  consider  the extended formulation \eqref{minlp2} of a general NLP problem and focus on the associated lifted set $\cSl$ in \eqref{eq.lift}.  We show a lifting result on constructing maximal $\cSl$-free sets.

Let $z\deq(x,y)$ denote the vector variable in the extended formulation \eqref{minlp2}, with its index set being $[n+\ell]$. Consequently, we have $z_{[n]} = x$ and $z_{[n+1:n+\ell]} = y$. Consider a closed subset $\cX$ of the domain $\bigcap_{i \in [\ell]} \dom(g_i)$ for $x$, and let $\cY$ be a closed subset of the domain $\bigtimes_{i \in [\ell]} \rang(g_i)$ for $y$. The ground set $\cG$ can, thus, be set as $\cX \times \cY$. Consequently, the lifted set $\cSl$ in \eqref{eq.lift} admits the form $\{(x,y) \in \cG: y = g(x)\}$.

Given that each $g_i(x)$ (for $i\in[\ell]$) may only depend on a subset of  variables indexed by $J_i \subseteq [n]$, we can express $g_i(x)$ as a lower order function $g'_i(x_{J_i})$ defined over $\bR^{J_i}$. Let $I_i \deq J_i \cup \{i + n\}$, and denote its complement by $I_i^c \deq [n+\ell] \smallsetminus I_i$. As above, we consider a closed subset $\cX^i$ of $\dom(g'_i)$ and $\cY^i$ of $\rang(g'_i)$. Consequently, the graph, epigraph, and hypograph of $g'_i$ reside within sets $\cG^i \deq \cX^i \times \cY^i$, \eg $\epi(g'_i)=\{(x_{J_i}, y_i) \in \cG^i: g'_i(x_{J_i}) \le y_i\}$.

We refer to $\cX,\cY, \{\cX^i, \cY^i\}_{i \in [\ell]}$ as the \emph{underlying sets} of the lifted set $\cSl$. The sets are said to be \emph{1d-convex decomposable} by a collection $\{\cD_j \}_{j \in [n+\ell]}$ of closed  convex sets in $\bR$, if  $\cX= \bigtimes_{j \in [n]}\cD_j, \cY= \bigtimes_{j \in [n+1:n+\ell]}\cD_j$, and, for all $i \in [\ell]$, $\cX^i = \bigtimes_{j \in J_i}\cD_j, \cY^i = \cD_{n+i}$. This decomposability condition restricts the domains to Cartesian products of real lines, intervals, or half lines, thereby excluding complicated domain structures.

The decomposability condition allows the analysis of sets with fewer variables. The construction of $\epi(g'_i)$-free sets and $\hyp(g'_i)$-free sets is in general simpler than the construction of $\cSl$-free sets. We show that any maximal $\epi(g'_i)$-free or $\hyp(g'_i)$-free set can be transformed into a maximal $\cSl$-free set.
\begin{theorem}
\label{thm.decomp}
Suppose the underlying sets of $\cSl$ are 1d-convex decomposable and $g$ is continuous. For some $i \in [\ell]$, let $\cC$ be a maximal $\epi(g'_i)$-free set or a maximal $\hyp(g'_i)$-free set in $\cG^i$. Then, the lifted set $\bar{\cC}\deq \cC \times \bR^{I^c_i}$ is a maximal $\cSl$-free set in $\cG$, where $\bR^{I^c_i}$ is the $|I^c_i|$-dimensional Euclidean space indexed by $I^c_i$.
\end{theorem}

See the proof in the appendix. For any $i\in[\ell]$, we call the operation $\cC \times \bR^{I^c_i}$ the \emph{orthogonal lifting} of $\cC$ with respect to $g_i$. A similar lifting result for integer programming is given by \cite[Lemma~4.1]{conforti2011}: given $\cS \deq \bZ^{n_1} \times \bR^{n_2}$, any maximal lattice-free set (\ie $\bZ^{n_1}$-free set) can be transformed into a maximal $\cS$-free set by orthogonal lifting. Therefore, \Cref{thm.decomp} serves as the NLP counterpart to this lemma (whose proof is also similar). This theorem allows us to focus on low-dimensional projections of the lifted set.
 We will show in \Cref{cor.maxsig} that the signomial lift satisfies the prerequisites of  \Cref{thm.decomp}.
The following example illustrates the application of \Cref{thm.decomp}.

\begin{example}
 Consider a lifted set $\cSl$ defined as $$\{(x_1,x_2,x_3,x_4, y_1,y_2, y_3): y_1= \exp(x_1 - x_2 /x_3)\land  y_2 = \log(x_1) \land y_3 = \sin(x_1 / x_4)\}.$$

One can verify that the 1d-convex decomposable condition holds for $\cD_1 = \bR_+$,  $\cD_j = \bR$ (for $j \in [2:7]$). Then $\cG\deq \bR^1_+ \times \bR^6$.  We use $\log(x_1)$ to construct a $\cSl$-free set. A maximal $\cSl$-free set can be $\{(x_1,x_2,x_3,x_4, y_1,y_2, y_3) \in \cG:  y_2 \le \log(x_1)\}$. Since $\log(x_1)$ is defined over positive  reals, this example gives a reason to restrict maximality over $\cG$.
\end{example}

\subsection{Sufficient conditions on maximality}

We provide sufficient conditions for the maximality of $\cS$-free sets for two general classes of nonconvex sets $\cS$. At the beginning, we give an overview of some basic results of convex analysis. Our subsequent exposition relies on the use of support functions of convex sets. The properties of support functions can be summarized as follows.
\begin{lemma}\cite[Chap. C]{hiriart2004fundamentals}
\label{lem.supconv}
 For  a full-dimensional closed convex set  \(\cC \subsetneq \bR^p\), let \(\sigma_{\cC}: \bR^p \to \bR, \lambda \mapsto \sup_{z \in \cC}{\lambda}\cdot z\) be the support function of \(\cC\). Then: (i) \label{lem.supconv2} \(\cC = \{z \in \bR^p:\forall {\lambda} \in \dom(\sigma_{\cC}) \textup{ } {\lambda}\cdot z\le \sigma_{\cC}({\lambda}) \}\), (ii)  \(\inter(\cC) = \{z \in \bR^p: \forall {\lambda} \in  \dom(\sigma_{\cC})  \smallsetminus\{0\} \textup{ } {\lambda}\cdot z <\sigma_{\cC}({\lambda}) \}\), (iii)
	$\sigma_{\cC}(\rho \lambda) = \rho \sigma_{\cC}(\lambda)$ for any $\rho > 0$. Moreover, for any closed convex set \(\cC'\) including $\cC$,   $\sigma_{\cC} \le \sigma_{\cC'}$.
\end{lemma}

 A valid inequality $a \cdot z \le b$ of $\cC$ is  called a \emph{supported valid inequality}, if there exists a \emph{supporting point} $z' \in \bd(C)$ such that \(a \cdot z' =  b\).  Geometrically, a closed convex set is the intersection of half-spaces associated with supported valid inequalities.

 \begin{observation}
  It follows from \Cref{lem.supconv} that every supported valid inequality of $\cC$ must admit the form $\lambda \cdot z\le \sigma_{\cC}(\lambda)$  for some $\lambda \in \dom(\sigma_{\cC})$, where the supremum  $\sigma_{\cC}(\lambda)$ is attained at its supporting points.
 \end{observation}

An inequality of the form $\lambda \cdot z \leq \sigma_{\cC}(\lambda)$, for $\lambda \in \dom(\sigma_{\cC})$, is referred to as an \emph{exposed valid inequality}, if there exists an \emph{exposing point} $z' \in \bd(C)$ such that $\lambda \cdot z' = \sigma_{\cC}(\lambda)$ and,  for all \({\lambda'} \in \dom(\sigma_\cC) \smallsetminus\{{\rho\lambda}\}_{\rho > 0}\), \({\lambda'}\cdot z' < \sigma_{\cC}(\lambda')\).
 
  \begin{observation}
  \label{obs.1}
An exposed valid inequality must be a supported valid inequality. Conversely, a supported valid inequality is an exposed valid inequality if the manifold $\bd(\cC)$ is smooth at its supporting point. For example, $\cC_1\deq\{(x,y) \in \bR^2: y = x^2\}$ is a smooth manifold, so any supported valid inequality of $\cC_1$ is exposed; $\cC_2\deq\{(x,y) \in \bR^2: y = |x|\}$ is smooth at $x \in [1,2]$, so any supported valid inequality of $\cC_2$ with support point $(x,y)$ ($x \in [1,2]$) is also exposed by the same point; however, a supported valid inequality of $\cC_2$ with supporting point $(x,y)$ ($x = 0$) cannot be exposed, since there are infinitely many supported valid inequalities at the same point.
 \end{observation}

The first lemma we present holds for full-dimensional nonconvex sets $\cS$. As shown in \Cref{fig.sdef1,fig.sdef2}, we have observed the geometric equivalence between the closed convex inner approximation of $\cl(\cS^c)$ and $\cS$-free sets. The  lemma provides a sufficient condition for the maximality of closed convex inner approximations.

\begin{lemma} \cite[Adapted from Thm.~3.1]{munoz2020}\label{prop.max}
	Let \(\cF\) be a full-dimensional closed set in \(\bR^p\), and let \(\cC  \subseteq  \cF\) be a full-dimensional closed convex set. If, for any \(z^\ast \in \inter(\cF \smallsetminus \cC)\) and any  $\lambda \in \dom(\sigma_{\cC})$ such that $\lambda \cdot z^\ast > \sigma_{\cC}(\lambda)$, there exists a point $z' \in \bd(\cF) \cap \bd(\cC)$ exposing $\lambda \cdot z \le \sigma_{\cC}(\lambda)$, then $\cC$ is a maximal convex inner approximation of $\cF$.
	\end{lemma}

We call  $z^\ast$ in \Cref{prop.max} an \emph{outlier point}, by which we try to enlarge an $\cS$-free set, and let the scope $L(z^\ast):=\{\lambda \in \dom(\sigma_{\cC}):\lambda \cdot z^\ast > \sigma_{\cC}(\lambda)\}$ identify the strictly separating valid inequalities for $z^\ast$.  Thm.~3.1 in \cite{munoz2020} has a different quantification than \Cref{prop.max}: it does not quantify $z^\ast$, and it requires the scope of $\lambda$ to be a subset $\Gamma \subseteq \dom(\sigma_{\cC})$, which is declared according to the context. Instead, \Cref{prop.max}  quantifies $\lambda$ explicitly, whose scope $L(z^\ast)$ depends on $z^\ast$. Thus,  \Cref{prop.max} allows, for each point $z^\ast$, having different scope $L(z^\ast)$ of $\lambda$.  One can prove \Cref{prop.max} by adapting the proof for \cite[Thm.~3.1]{munoz2020}. For self-completeness, we give a proof in the appendix using support functions.

We next focus on a specific type of function, namely  \emph{positive homogeneous functions}.  We summarize their properties as follows.
\begin{lemma}
	\label{lem.homo}
Let  $f$  be a positive homogeneous function of degree $d \in \bR$, such that, for any $z \in \dom(f) \subseteq \bR^p$ and any $\rho \in \bR_{++}$, $f(\rho z)=\rho^d f(z)$. Then: (i) $\inter(\dom(f))$ is a cone, and (ii) if $d=1$, then for any $\breve{z} \in \dom(f)$, \(\lin{f}{\breve{z}}(z) =\nabla  {f}(\breve{z})\cdot z\) for $z \in \dom(f)$ and \(\lin{f}{\breve{z}}(z)  = f(z)\) for \( z = \rho \breve{z}\) with \(\rho \in \bR_{++}\).
\end{lemma}

The proof is in the appendix. We recall that $\lin{f}{\breve{z}}$ in the above lemma is defined in Eq.~\eqref{def.Xi}. Moreover, $\dom(f)$ is embedded in $\bR^p$, so we call $\bR^p$ the \emph{ambient space} of $f$.

The second theorem we present offers a more structured result, specifically related to nonconvex DCC sets $\cS$. \cite[Thm.~5.48]{serranothesis} provides a sufficient condition for the maximality of the $\cS$-free set described in \Cref{cor.dc}. However, to clearly distinguish it from our result below, we translate the condition into our setting as follows: (i) the functions $f_1$ and $f_2$ are superlinear, i.e. they are positive homogeneous of degree 1 and superadditive (note that superlinear functions are concave), (ii) they are separable and act independently on different variables $u$ and $v$, (iii) $f_1$ is negative everywhere except at 0, (iv) the linearization point $\relx{v}$ of $f_2$ is nonzero, and (v) the domains $\dom(f_1)$ and $\dom(f_2)$ are Euclidean spaces.

Our second theorem provides an alternative condition for maximality that relaxes condition (i) by requiring only that one of $f_1$ or $f_2$ be positive homogeneous of degree 1, while imposing mild regularity conditions. Moreover, the domains can be full-dimensional convex cones.

\begin{theorem}
\label{thm.max}
	For every $i \in \{1,2\}$, let  $f_i$ be   concave.  Let $\cS \deq \{(u,v) \in \dom(f_1) \times \dom(f_2): f_1(u) - f_2(v) \le 0\}$. Suppose that: (i) at least one  of $f_1,f_2$ is positive homogeneous of degree 1, (ii) $f_1,f_2$ are both positive/negative over the interiors of their domains, (iii) $f_1$ is continuously differentiable over $\inter(\dom(f_1))$, and  (iv) $\dom(f_1),\dom(f_2)$ are full-dimensional in the ambient spaces of $f_1,f_2$, respectively. Then, for any $\relx{v} \in \inter(\dom(f_2))$, $\cC \deq \{(u,v) \in \dom(f_1) \times \dom(f_2): f_1(u) -\lin{f_2}{\relx{v}}(v) \ge 0\}$ is maximally $\cS$-free in $\dom(f_1) \times \dom(f_2)$.    
\end{theorem}
\begin{proof}
	We first adapt \Cref{cor.dc} by restricting the domain of $z$ to the convex ground set $\cG\deq\dom(f_1) \times \dom(f_2)$. It follows from \Cref{cor.dc} that $\cC$ is  an $\cS$-free set in $\cG$. Since $\dom(f_1) \times \dom(f_2)$ are full-dimensional, $\cS,\cC, \cG$ are  full-dimensional. As $\cS,\cC \subseteq \cG$, the maximality of $\cC$ in $\cG$ is equivalent to that $\cC$ is a maximal convex inner approximation of $\cF\deq\cl(\cS^c) \cap \cG = \{(u,v) \in \cG: f_1(u) -f_2(v) \ge 0\}$. Note that $\cF$ is  full-dimensional. We then apply \Cref{prop.max} to prove that $\cC$ is a maximal convex inner approximation of $\cF$.  Let \(z^\ast \in \inter(\cF \smallsetminus \cC)\) be any outlier point. It follows from the separating hyperplane theorem that there exists a supported valid inequality $\lambda \cdot z \le \sigma_{\cC}(\lambda)$  of $\cC$ such that $\lambda \cdot z^\ast > \sigma_{\cC}(\lambda)$. Since $\cF \smallsetminus \cC \subseteq \cG$, $ \inter(\cF \smallsetminus \cC) \subseteq \cG$. Since $\cC \subseteq \cG$, the inequality cannot be supported by a valid inequality at $\bd(\cG)$, so the inequality must be a valid inequality supported at $\cC \smallsetminus \bd(\cG)$. It follows from the concavity of $f_1$ that the inequality must admit the form $\lin{f_1}{\breve{u}}(u) -\lin{f_2}{\relx{v}}(v) \ge 0$ for some $\breve{u} \in \dom(f_1)$ (identical up to a positive multiplier). By the smoothness of $f_1$, w.l.o.g, we can perturb $\breve{u}$ such that it is in $ \inter(\dom(f_1))$. Let $\breve{v} := \relx{v}$. We now have that $\breve{u} \in \inter(\dom(f_1)), \breve{v} \in \inter(\dom(f_2))$.  We will prove that $\lin{f_1}{\breve{u}}(u) -\lin{f_2}{\breve{v}}(v) \ge 0$ is exposed by a point $(u',v') \in (\bd(\cF) \cap \bd(\cC)) \cap \inter(\cG)$. It suffices to show that the following three equations hold:
	\begin{equation}
	\label{eq.three}
    	\begin{aligned}
        	\lin{f_1}{\breve{u}}(u') -\lin{f_2}{\breve{v}}(v') &= 0 \quad \textup{(\ie supported at $(u',v')$)},\\
        	f_1(u') -\lin{f_2}{\breve{v}}(v') &= 0  \quad \textup{(\ie $(u',v') \in \cC$)},\\
        	f_1(u') -f_2(v') &= 0 \quad  \textup{(\ie $(u',v') \in \cF$)}.
    	\end{aligned}
	\end{equation}
   Since $\cC \subseteq \cF$ and they  are both full-dimensional, the last two equations imply that $(u',v') \in \bd(\cC) \cap \bd(\cF)$. As $f_1$ is continuously differentiable and concave in the interior of its domain, the graph of $f_1(u) -\lin{f_2}{\breve{v}}(v)$  over $\inter(\cG)$ is a smooth manifold embedded in $\inter(\cG) \times \bR$. The intersection of a smooth manifold with a hyperplane yields another lower-dimensional smooth manifold. This implies that the  level set $\cC$ of $f_1(u) -\lin{f_2}{\breve{v}}(v)$ is also smooth at any point  $(u,v) \in \inter(\cG) \cap \cC$. By \Cref{obs.1}, $(u,v)$  is an exposing point.  Since $(u',v') \in \cC \cap  \inter(\cG)$, $(u',v')$ is an exposing point, and the maximality of $\cC$ is verified. We now proceed to construct $(u',v')$ from $(\breve{u},\breve{v})$ and prove \eqref{eq.three}.
  Let $\rho \deq  f_2(\breve{v}) / f_1(\breve{u})$.  Since $\breve{u} \in \inter(\dom(f_1)), \breve{v} \in \inter(\dom(f_2))$, by the assumption, $\rho  > 0 $.  We consider the following two cases separately.
  
  \textbf{Case i.} We first suppose that $f_1$ is positive homogeneous of degree 1.  Let $(u', v')\deq (\rho \breve{u}, \breve{v})$, which,  by \Cref{lem.homo}, is in $\inter(\cG)$.
We have that:
$$ 
  	f_1(u') \stackrel{(i.1)} = \lin{f_1}{\breve{u}}(u')
    	\stackrel{(i.2)}= \rho f_1(\breve{u})
    	\stackrel{(i.3)}= f_2(\breve{v})
    	\stackrel{(i.4)}= f_2(v') \stackrel{(i.5)}= \lin{f_2}{\breve{v}}(v'),
 $$ 
where equations $(i.1),(i.2)$ follow from \Cref{lem.homo},   $(i.3)$ follows from the definition of $\rho$, and $(i.4),(i.5)$  follow from $v'=\breve{v}$.
 
 \textbf{Case ii.}
We then suppose that $f_2$ is positive homogeneous of degree 1.  Let $(u', v')\deq ( \breve{u}, \breve{v} / \rho) \in \inter(\cG)$.
We have that:
$$
\lin{f_1}{\breve{u}}(u')   \stackrel{(ii.1)} =   f_1(u') \stackrel{(ii.2)} =  f_1(\breve{u})  
    	 \stackrel{(ii.3)} =  f_2(\breve{v}) /  \rho      	 \stackrel{(ii.4)} =   f_2(v') 
    	 \stackrel{(ii.5)} = \lin{f_2}{\breve{v}}(v'),
$$
	where equations $(ii.1), (ii.2)$ follow from $\breve{u} = u'$,  $(ii.3)$ follows from the definition of $\rho$, and  $(ii.4),(ii.5)$ follow  from \Cref{lem.homo}.  Therefore, \eqref{eq.three} are satisfied in both cases.
\end{proof}

We present the motivation for restricting the maximality of the set $\cC$ within the ground set $\dom(f_1) \times \dom(f_2)$. The main reason for this restriction arises from the difficulty of finding a nontrivial concave extension of $f_1$ over its ambient space such that for all $u \notin \dom(f_1)$, $f_1(u) > -\infty$. While such an extension can exist geometrically, the construction of a closed expression remains unclear. In the next section, we will examine a specific example to illustrate this point.

Moreover, we will apply the above theorem to develop DCC formulations for a nonconvex set. In particular, the functions $f_1$ and $f_2$ must not simultaneously have positive homogeneity of degree 1, and their domains are non-negative orthants. Consequently, the relaxed condition for homogeneous degrees and domains in \Cref{thm.max} becomes necessary. We give two examples for verification \Cref{thm.max}.
\begin{example}
    Let $f_1(u) \deq u$ with $\dom(f_1) \in \bR$, and let $f_2(v) \deq \sum_{i \in [n]}\sqrt{v_i}$ with $\dom(f_2) = \bR^n_{+}$. Note that $f_1,f_2$ are concave, $\dom(f_2)$ is a non-negative orthant, and $f_1$ is positive homogeneous of degree 1. Let $\cG \deq \bR \times \bR^n_{+}$. One can verify that the presupposition of \Cref{thm.max}  is satisfied. Then, $\cS \deq \{(u,v) \in \cG: u -   \sum_{i \in [n]}\sqrt{v_i}\le 0\}$ is a convex set. It is easy to see that $\cC \deq \{(u,v) \in \cG: u -  \sum_{i \in [n]}(\sqrt{\relx{v}_i} + (v_i - \relx{v}_i)/ \sqrt{\relx{v}_i}) \ge 0\}$ is maximally $\cS$-free in $\cG$ with $\relx{v} > 0$.
\end{example}
\begin{example}
Exchange the functions $f_1,f_2$ in the previous examples. Then, $\cS \deq \{(u,v) \in \cG:    \sum_{i \in [n]}\sqrt{v_i} - u\le 0\}$ is a reverse-convex set. It is easy to see that $\cC \deq \{(u,v) \in \cG: \sum_{i \in [n]}\sqrt{v_i} - u \ge 0\}$ is the unique maximal $\cS$-free  set in $\cG$.
\end{example}

\section{Signomial-lift-free sets and intersection cuts}
\label{sec.icforsp}
In this section, we  construct  (maximal) signomial-lift-free sets and generate intersection cuts for SP.  

\subsection{Signomial-lift-free and signomial-term-free sets}
\label{sec.siliftfree}

We introduce and study new formulations of signomial-term sets. We transform signomial-term sets into DCC sets. We also construct signomial term-free sets and lift them to signomial term-lift-free sets. The maximality of these sets is studied, and a comparison is made between signomial term-free sets derived from different DCC formulations.

We consider an $n$-variate signomial term $\psi_{\alpha}(x)$ arising in the extended formulation \eqref{minlp2} of SP.
The exponent vector $\alpha$ may contain negative/zero/positive entries. We extract  two sub-vectors $\alpha_{-}$ and $\alpha_{+}$ from $\alpha$ such that $\alpha_{-}  \in \bR^{\eta}_{--}$ ($\eta$-dimensional negative orthant) and $\alpha_{+} \in \bR^{\kappa}_{++}$ ($\kappa$-dimensional positive orthant), and let $x_- \in \bR^{\eta}$ and  $x_+ \in \bR^{\kappa}$ be the corresponding sub-vectors of $x$. Entries $x_j$ with $\alpha_j = 0$ are excluded from consideration, and so $\eta + \kappa$ may be smaller than $n$. Since $\psi_{\alpha}(x)$ only depends on $x_-$ and $x_+$, it can be represented in the form of $x_{-}^{\alpha_{-}} x_{+}^{\alpha_{+}}$ of lower order.

Let $\lessgtr$ (resp. $\lesseqgtr$) denote  $<$ or $>$ ($\le$ or $\ge$). 
We consider the \emph{signomial-term set} as epigraph or hypograph of $x_{-}^{\alpha_{-}}  x_{+}^{\alpha_{+}}$:
\begin{equation}
	\label{eq.sigset}
  \cSt = \{(x_-,x_+, t) \in \bR_{+}^{\eta + \kappa+1}: t  \lesseqgtr x_{-}^{\alpha_{-}}  x_{+}^{\alpha_{+}} \}.
\end{equation}

We first give DCC reformulations of  signomial-term sets.
The interior of $\cSt$ in \eqref{eq.sigset} is
 \begin{equation*}
   \inter(\cSt) = \{(x_-,x_+, t) \in \bR_{++}^{\eta + \kappa+1} : t \lessgtr x_{-}^{\alpha_{-}}  x_{+}^{\alpha_{+}} \}.
\end{equation*}
Reorganizing the signomial terms and taking the closure of the set, we recover
\begin{equation*}
  	\cSt =\{ (x_-,x_+, t) \in \bR_{+}^{\eta + \kappa+1}:  t x_-^{-\alpha_{-}} \lesseqgtr x_+^{\alpha_{+}}\}.
\end{equation*}

Notably, the exponents associated with signomial terms on both sides are now strictly positive. Let $u \deq (t,x_-), v \deq x_+$, let $h \deq \eta+1$, and let $k \deq \kappa$. Then, $\psi_{\beta'}(u)=t x_-^{-\alpha_{-}}$ and $ \psi_{\gamma'}(v)=x_+^{\alpha_{+}}$, where \(\beta'\deq(1, -\alpha_{-}) \in \bR_{++}^h\) and \(\gamma' \deq \alpha_{+}  \in \bR_{++}^k\).  After the change of variables, the set  admits the following form:
\begin{equation}
\label{eq.form1}
	\cSt =  \{(u, v) \in \bR_{+}^{h + k}:\, \psi_{\beta'}(u) \lesseqgtr \psi_{\gamma'}(v)\}.
\end{equation}

The formulation \eqref{eq.form1} exhibits symmetry between $u$ and $v$.  We can therefore consider w.l.o.g.~the inequality ``$\le$'' throughout the subsequent analysis.
Since the signomial terms $\psi_{\beta'}(u), \psi_{\gamma'}(v)$
are non-negative over $\bR_{+}^h, \bR_{+}^k$, we can take any positive power $\mu \in \bR_{++}$ on both sides of \eqref{eq.form1}.  Finally, the signomial  term set in \eqref{eq.sigset} admits the following form:
\begin{equation}
\label{eq.form2}
	\cSt =  \{(u, v) \in \bR_{+}^{h + k}:\, \psi_{\beta}(u) - \psi_{\gamma}(v) \le 0\},
\end{equation}
where $\beta \deq \mu \beta'$, and $\gamma \deq \mu \gamma'$.

A signomial term $\psi_{\alpha}(x)$ is said to be a \emph{power function} if $\alpha \ge 0$, and $\lVert \alpha \rVert_1 \le 1$. According to \cite{aps2018mosek,chares2009cones}, power functions are concave over the non-negative orthant; if additionally $\lVert \alpha \rVert_1 = 1$, $\psi_{\alpha}(x)$ is positive  homogeneous of degree 1. Moreover, $\psi_{\alpha}(x)$  has an extended exponential cone representation \cite{exponentialcpne}, which gives another proof of its convexity. Through an appropriate scaling of the parameter $\mu$, we obtain a family of DCC reformulations \eqref{eq.form2} of signomial-term sets. We let $\cG\deq  \bR_{+}^{h + k}$, and use the reverse-minorization technique  to construct signomial-term-free sets. We recall that the definition of the operator $\Xi$ is given in Eq.~\eqref{def.Xi}.
\begin{proposition}
\label{prop.maxsig}
Let $\max(\lVert \beta \rVert_1, \lVert \gamma \rVert_1 ) \le  1$. For any $\relx{v} \in \bR^{k}_{++}$, \begin{equation}
\label{eq.sigfree}
	\cC \deq  \{(u, v) \in \bR_{+}^h \times \bR^{k}:\, \psi_{\beta}(u) - \lin{\psi_{\gamma}}{\relx{v}}(v) \ge 0\}
\end{equation} is a signomial-term-free ($\cSt$-free) set.   If $\max(\lVert \beta \rVert_1, \lVert \gamma \rVert_1 ) = 1$, then $\cC$ is a maximal signomial-term-free set in $\cG$.
 \end{proposition}
 \begin{proof}
Since $ \max(\lVert \beta \rVert_1, \lVert \gamma \rVert_1 ) \le 1 $,  $ \psi_{\beta}(u), \psi_{\gamma}(v)$ are concave. By \Cref{cor.dc}, $\cC$ is signomial-term-free. If $\max(\lVert \beta \rVert_1, \lVert \gamma \rVert_1 ) = 1$, then  at least one of $\lVert \beta \rVert_1,  \lVert \gamma \rVert_1$ is 1. Therefore, one of $\psi_{\beta}(u), \psi_{\gamma}( v)$ is positive  homogeneous of degree 1.  Moreover, $\psi_{\beta}(u), \psi_{\gamma}( v)$ are both continuously differentiable and positive over positive orthants $\bR^h_{++},\bR^{k}_{++}$ (the interiors of their domains). Since $\cG = \dom(\psi_{\beta}) \times  \dom(\psi_{\gamma})$,  by \Cref{thm.max}, $\cC \cap \cG= \{(u, v) \in \cG:\, \psi_{\beta}(u) - \lin{\psi_{\gamma}}{\relx{v}}(v) \ge 0\}$  is  a maximal signomial-term-free set in $\cG$. Therefore, $\cC$ is also a maximal signomial-term-free set in $\cG$.
 \end{proof}
 
Given that $ \max(\lVert \beta \rVert_1, \lVert \gamma \rVert_1 ) = 1$ results in a desirable DCC formulation for the signomial-term set, we refer to this formulation as its \emph{normalized DCC formulation}.  Comparing \Cref{prop.maxsig} to \Cref{thm.max}, we extend the domain of $\lin{\psi_{\gamma}}{\relx{v}}(v)$ from $\bR^{k}_{+}$ to $\bR^{k}$, since it is an affine function. However, the further extension requires a non-trivial concave extension of the power function $\psi_{\beta}$, which we are unaware of.

We have reduced the $n$-variate signomial term $\psi_{\alpha}(x)$ to a signomial term $ x_{-}^{\alpha_{-}}  x_{+}^{\alpha_{+}}$ of lower order and constructed the corresponding signomial-term-free sets. A similar reduction is observed for $g_i$ to $g'_i$ in \Cref{sec.lift}, where we demonstrate the relationship between  $\cSl$-free sets and  $\epi(g'_i)$-free/$\hyp(g'_i)$-free sets.

Next, we let the lifted set $\cSl$ be the signomial lift,  where all $g_i$ are signomial terms. Each equality constraint $y_i = g_i(x)$ defining the signomial lift is equivalent to two inequality constraints $y_i \lesseqgtr g_i(x)$.  Applying the normalized DCC reformulation to  these inequality constraints, we thus obtain a reformulation of the signomial lift, which we call its \emph{normalized DCC reformulation}.

\begin{corollary}
\label{cor.maxsig}
 Let $\cC$ be as in \eqref{eq.sigfree}, where  $\psi_{\alpha} = g_i$ for some $i \in [\ell]$ and $\max(\lVert \beta \rVert_1, \lVert \gamma \rVert_1 ) = 1$.  Then the orthogonal lifting of $\cC$ w.r.t. $g_i$ is a maximal signomial-lift-free ($\cSl$-free) set in the non-negative orthant.
\end{corollary}
\begin{proof}
	We verify that the conditions of \Cref{thm.decomp} are satisfied by the signomial lift. For any $i \in [\ell]$, the signomial term $g_i$ is continuous, and its domain and range are  $\bR_{++}$. Let $J_i$ be the index set of variables of its reduced signomial term $g'_i$. Let $\cX \deq \bigtimes_{j \in [n]} \bR_{++}, \cY \deq  \bigtimes_{j \in [\ell]} \bR_{++}$. For all $j \in [n+\ell]$,  let $\cD_j \deq \bR_{++}$. For all $i \in [\ell]$, let $\cX^i \deq \bigtimes_{j \in J_i} \bR_{++}, \cY^i \deq \bR_{++}$. The underlying sets of the signomial lift are $\cX,\cY, \{\cX^i, \cY^i\}_{i \in [\ell]}$ that are 1d-convex decomposable by $\{\cD_j\}_{j \in [n+\ell]}$.  By \Cref{prop.maxsig}, $\cC$  is a maximal $\hyp(g'_i)$-free set in $\cX^i \times \cY^i$. By \Cref{thm.decomp}, its orthogonal lifting w.r.t. $g_i$ is a maximal signomial-lift-free set in positive orthant. By continuity of $\psi_{\beta}, \psi_{\gamma}$, we  change the ground set (the positive orthant) to its closure, \ie non-negative orthant.
\end{proof}

The following examples show signomial term-free sets from different DCC formulations.

\begin{figure}[!ht]
 \centering
  \subfloat[$\cSt$ and $\cC_{1}$.]{\includegraphics[width=0.3\textwidth]{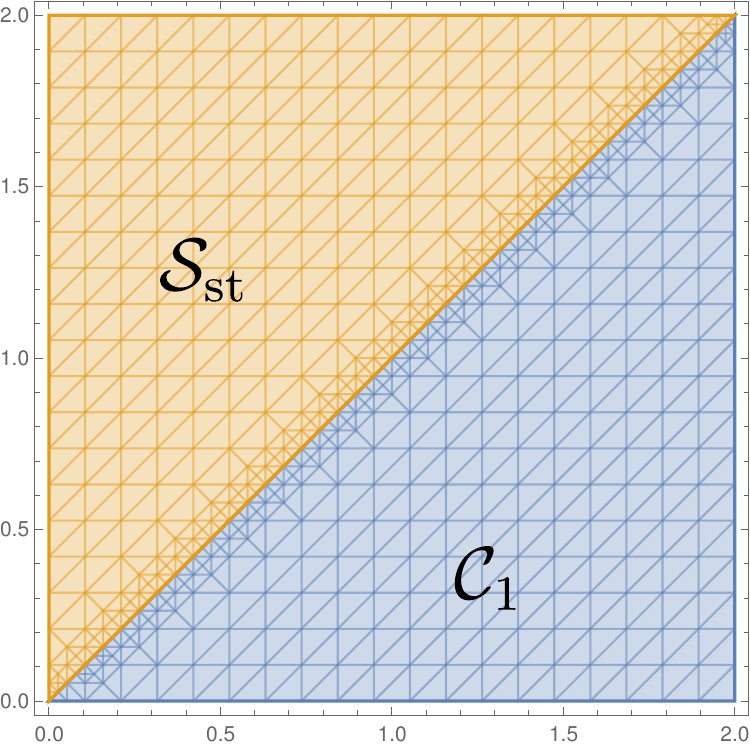}}\label{fig.ideal.log}
  \hfill
  \subfloat[$\cSt$ and $\cC_{2}$.]{\includegraphics[width=0.3\textwidth]{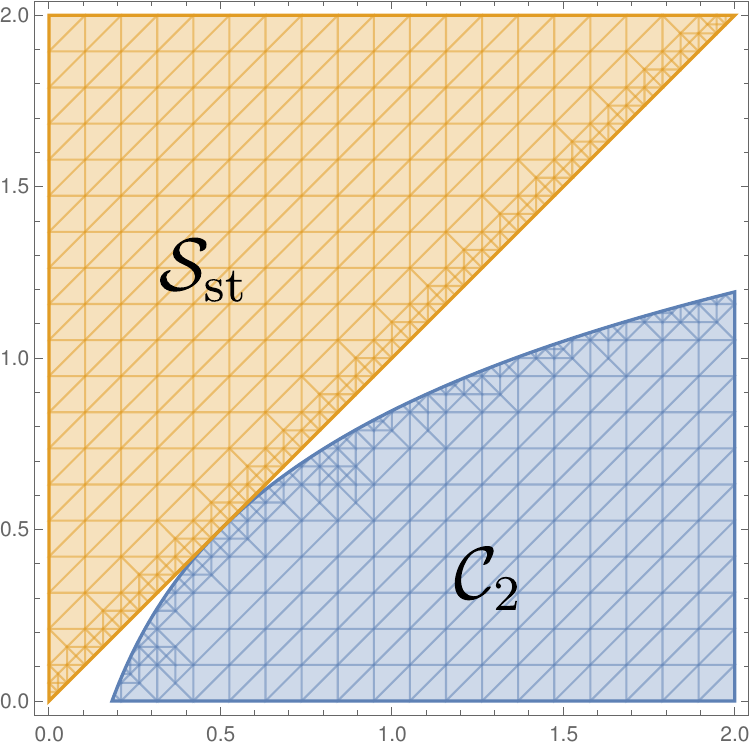}}\label{fig.ideal.power}
\hfill
  \subfloat[$\cSt$ and $\cC_{3}$.]{\includegraphics[width=0.3\textwidth]{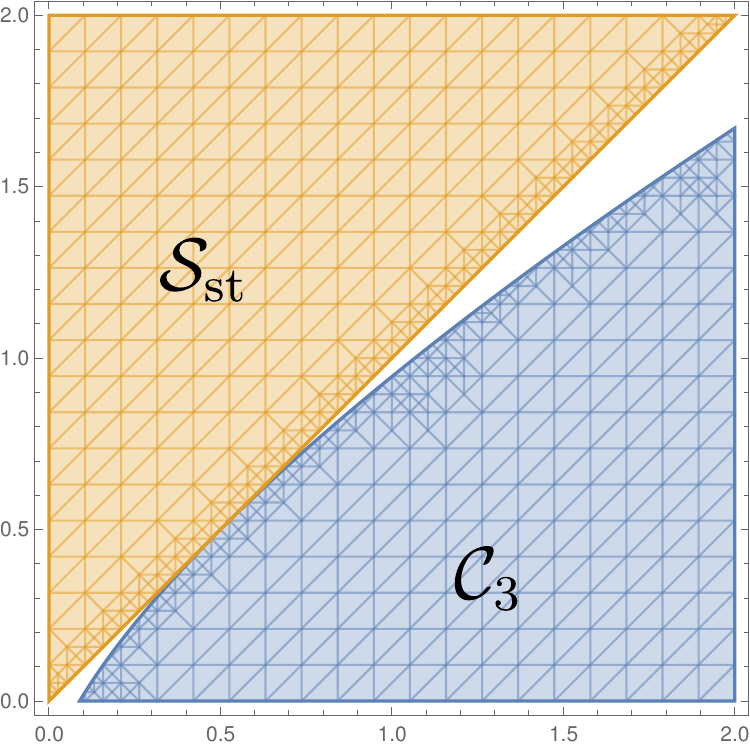}}\label{fig.exmpl.linear}
  \caption{$\cSt$-free sets from \Cref{sec.idealandnonideal}.}
  \label{fig.ideal}
\end{figure}

\begin{example}[Comparison of DCC formulations]
\label{sec.idealandnonideal}
Consider $
   \cSt =\{(u, v) \in \bR^2_{+}: u \le v\}$, which is already in normalized DCC formulation.
It is easy to see that $
   \cC_1 \deq\{(u, v) \in \bR_{+} \times \bR: u \ge v\}$ is a maximal $ \cSt$-free set in $\bR^2_{+}$ given by \Cref{prop.maxsig}. Let $\relx{v} \in \bR_{++}$ be a linearization point. Consider the set  $\cSt' \deq \{(u, v) \in \bR^2_{++}: \log(u) \le \log(v)\}$. We find that $\cSt' \subsetneq \cSt$, but two sets almost coincide except for some boundary points of $\cSt$. Since $\cSt'$  admits a DCC formulation, applying the reverse-minorization technique at $\relx{v}$ yields $ \cC_2 \deq\{(u, v) \in \bR^2_{+}: \log(u) - (\log(\relx{v}) + (v - \relx{v}) / \relx{v}) \ge 0\}$, which is also an $ \cSt$-free set. For any $0 < \mu < 1$, $\cSt =\{(u, v) \in \bR^2_{+}: u^{\mu} \le  v^\mu\}$ is a DCC set, applying the reverse-minorization technique at $\relx{v}$ yields $ \cC_3 \deq\{(u, v) \in \bR^2_{+}: u^{\mu} - ((1- \mu)\relx{v}^{\mu} + \mu \relx{v}^{\mu - 1}v)  \ge 0\}$, which is also an  $ \cSt$-free set. However, $\cC_2,\cC_3$ cannot be maximal in $\bR^2_{+}$, because their intersections with $\bR^2_{+}$ are not  polyhedral. These sets are visualized in \Cref{fig.ideal} with a linearization point $\relx{v} = 0.5$ and scaling parameter $\mu= 0.7$.
\end{example}

\begin{figure}[!ht]
 \centering
  \begin{subfigure}[b]{0.35\textwidth}
     	\centering
     	\includegraphics[width=\textwidth]{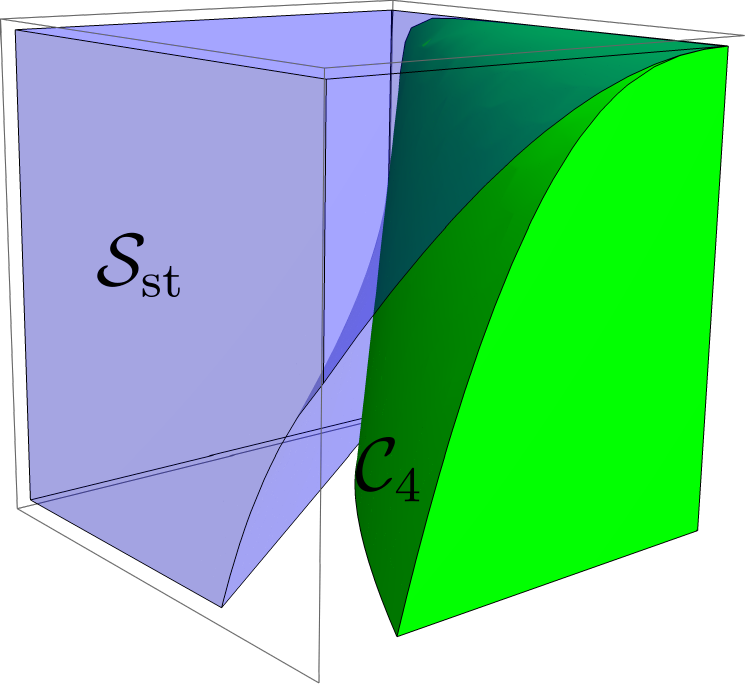}
     	\caption{$\cSt$ and $\cC_4$ from \Cref{example.1}.}
    	\label{fig.free1}
 \end{subfigure}
 \hfill
  \begin{subfigure}[b]{0.35\textwidth}
     	\centering
     	\includegraphics[width=\textwidth]{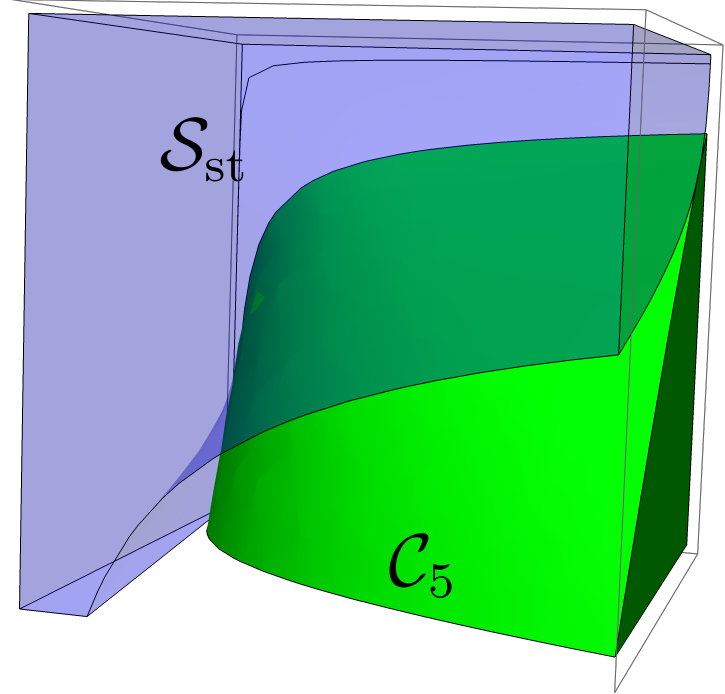}
     	\caption{$\cSt$ and $\cC_5$ from \Cref{example.2}.}
    	\label{fig.free2}
 \end{subfigure}
  \caption{$\cSt$ and $\cSt$-free sets from \Cref{example.1,example.2}.}
  \label{fig.free}
\end{figure}

\begin{example}
\label{example.1}
Consider the  hypograph of signomial term $x_1^{-2} x_2^{2}$ and
$
	\cSt = \{(x,y) \in \bR_{+}^3: y \le x_1^{-2} x_2^{2} \}.
$ For $(x,y) \in \bR_{++}^3$, $y \le x_1^{-2} x_2^{2}$ if and only if $ y^{1/3}x_1^{2/3} \le x_2^{2/3}$.  The following set is maximal $\cSt$-free in $\cG = \bR_+^3$:
$
	\cC_4 \deq \{(x,y) \in \bR_{+}^3:   y^{1/3}x_1^{2/3} \ge  \relx{x}_2^{2/3} + \frac{2}{3} \relx{x}_2^{-1/3} (x_2 - \relx{x}_2)\},
$
where $\relx{x}_2 \in \bR_{++}$. See \Cref{fig.free1} for $\relx{x}_2 = 0.2$.
\end{example}

\begin{example}
\label{example.2}
Consider the epigraph of signomial term $x_1^3 x_2$ and
$
	\cSt = \{(x,y) \in \bR_{+}^3: y \ge x_1^3 x_2\}.
$ For $(x,y) \in \bR_{++}^3$, $  y \ge x_1^3 x_2 $ if and only if $ y^{1/4} \ge x_1^{3/4}x_2^{1/4} $. The following set is maximal $\cSt$-free in $\cG = \bR_+^3$:
$
	\cC_5\deq \{(x,y) \in \bR_{+}^3:  \relx{y}^{1/4} + \frac{1}{4} \relx{y}^{-3/4}(y - \relx{y})\le x_1^{3/4}x_2^{1/4}  \},
$
where $\relx{y} \in \bR_{++}$. See \Cref{fig.free2} for $\relx{y} = 0.2$.
\end{example}

\subsection{Intersection cuts}
\label{sec.sepic}
We focus on the separation of intersection cuts for the extended formulation of SP. In \Cref{sec.prem} we presented a method to construct a simplicial cone $\cR$ from an LP relaxation. The vertex of this cone is a relaxation solution $\relx{z} = (\relx{x}, \relx{y})$. We choose $\relx{z}$ as the linearization point for applying the reverse-minorization technique.

We assume that the LP relaxation includes all linear constraints from \eqref{minlp2}. If $\relx{z}$ is infeasible for \eqref{minlp2}, then $\relx{z}$ does not belong to the signomial lift. Thus, there is a signomial term $g_i$ such that $\relx{y}_i \neq g_i(\relx{x})$. Given the reduced form $g'_i$, we obtain a set of signomial terms $\cSt$: If $g_i(\relx{x}) > \relx{y}_i$, we choose $\cSt$ to be the epigraph of $g'_i$; otherwise, we choose it to be the hypograph of $g'_i$. This signomial-term set yields a signomial term-free set $\cC$ in \eqref{eq.sigfree} containing $(\relx{u},\relx{v})$ in its interior (\Cref{cor.dc}). Using orthogonal lifting of \Cref{cor.maxsig}, we can transform $\cC$ into a signomial-lift-free set $\bar{\cC}$.

We next show how to construct an intersection cut in \eqref{eq.ic}. It suffices to compute  step lengths $\mu^\ast_j$ in \eqref{eq.iccoef} along extreme rays $r^j$ of $\cR$. Each step length $\mu^\ast_j$ corresponds to a boundary point $\relx{z}+ \mu^\ast_j r^j$ in $\bd(\bar{\cC})$. The left-hand-side $\psi_{\beta}(u) - \lin{\psi_{\gamma}}{\relx{v}}(v)$  of the inequality in \eqref{eq.sigfree} is a concave function over $(u,v) \in \bR_{+}^h \times \bR^{k}$. Its  restriction along the ray $\relx{z}+ \mu_j r^j \;(\mu_j \in \bR_+)$ is a univariate concave function:
\begin{equation*}
 \tau_j:\bR_+ \to \bR, \mu_j \mapsto  \tau_j(\mu_j) \deq \psi_{\beta}(\relx{u} +  r^j_u \mu_j)  - \lin{\psi_{\gamma}}{\relx{v}}(\relx{v}+ r^j_v\mu_j),
\end{equation*}
where $r^j_u$ and $r^j_v$ are the projections of $r^j$ on $u$ and $v$ respectively. Let
$	\bar{\mu}_j \deq \sup_{\mu_j \ge 0}\{\mu_j:\relx{u} +  r^j_u \mu_j \ge 0\}$.
 Therefore, $\mu^\ast_j$ is the first point in $[0,\bar{\mu}_j]$  satisfying the boundary condition: either $ \tau_j(\mu^\ast_j) = 0$ or $\mu^\ast_j = \bar{\mu}_j$. Since  $\tau_j$ is a univariate concave function and $\tau_j(0) > 0$, there is at most one positive point in $\bR_+$ where $\tau_j$ is zero. We employ the bisection search method \cite{press2007chapter} to find such $\mu^\ast_j$.

\section{Convex outer approximation}
\label{sec.outerforsp}

In this section we propose a convex nonlinear relaxation for the extended formulation \eqref{minlp2} of SP. This relaxation is easy to derive and allows us to generate valid linear inequalities, called outer approximation cuts, for SP. Unlike intersection cuts, outer approximation cuts do not require an LP relaxation \textit{a priori}, so solvers can employ them to generate an initial LP relaxation of \eqref{minlp2}.

With notation from \Cref{sec.siliftfree}, we additionally assume that  the domain of $u$ (resp. $v$) is a hypercube  $\cU$ (resp. $\cV$) in $\bR_{+}^{h}$(resp. $\bR_{+}^{k}$). The assumption fits with the common practice of MINLP solvers.  We construct the convex nonlinear relaxation by approximating each signomial-term set of the signomial lift within the hypercube. 

For brevity, we still call the intersection of  the set  in \eqref{eq.form2} and the hypercube $ \cU \times \cV$:
\begin{equation}
\label{eq.st2}
	\cSt \deq  \{(u, v) \in \cU \times \cV:\, \psi_{\beta}(u) - \psi_{\gamma}(v) \le 0\},
\end{equation}
 a signomial-term set.
 As long  as $ \max(\lVert \beta \rVert_1, \lVert \gamma \rVert_1 ) \le 1$, $\cSt$ is in a DCC formulation (in terms of the inequality constraint).

  We consider the normalized DCC formulation that has $\max(\lVert \beta \rVert_1, \lVert \gamma \rVert_1 ) = 1$. In \Cref{sec.alter,sec.conv}, we will explain the reason for choosing the  normalized DCC formulation. The signomial-term set is usually nonconvex, so our construction involves convexifying the concave function $\psi_{\beta}$ in \eqref{eq.st2}. This procedure yields a convex outer approximation of $\cSt$, which is non-polyhedral. Consequently, replacing $\cSt$ by its convex outer approximation, we obtain the convex nonlinear relaxation of \eqref{minlp2}.

Next, we introduce the procedure of relaxation. We should import the formal concepts of convex underestimators and convex envelopes.
Given a function $f$ and a closed set $\cD \subseteq \bR^p$, a convex function $f':\conv(\cD) \to \bR$ is called a convex underestimator of $f$ over $\cD$, if for all $x \in \cD$ $f'(x) \le f(x)$. The
 convex envelope  of $f$ is defined as the pointwise maximum convex underestimator of $f$ over $D$, and we denote it by $\conve_{\cD}(f)$.

In principle, the envelope construction procedure is similar to the convexification procedure of multilinear terms \cite{tawarmalani13}.
 The following lemma gives an extended formulation of the convex envelope of a concave function over a polytope, where the formulation is uniquely determined by the function values at the vertices of the polytope.
\begin{lemma}\cite[Thm.~3]{falk1976successive}
\label{thm.poly}
Let $P$ be a polytope in $\bR^n$, let $f:P \to \bR$ be a concave function over $P$, and let $Q$ be  vertices of $P$. Then,  $\conve_{P}(f)(x)=\min \{\sum_{q \in Q} \lambda_q f(q): \exists \lambda \in \bR^Q_{+},  \sum_{q \in Q} \lambda_q=1, \,  x =  \sum_{q \in Q} \lambda_q q\}$.
\end{lemma}

Based on the lemma above, we observe that the concave function $f$ is convex-extensible from its vertices (\ie $\conve_{P}(f)(x) = \conve_{Q}(f)(x)$ for $x \in P$), and $\conve_{P}(f)$ is a polyhedral function.

For the case of $P=\cU \deq \prod_{j \in [h]} [\underline{u}_j, \overline{u}_j]$ and $f =\psi_{\beta}$, $Q = \{q \in \bR^h: \forall j \in [h] \; q_j = \underline{u}_j \lor q_j = \overline{u}_j\}$ is the set of vertices of the hypercube $\cU$.  The lemma yields an extended formulation of $\conve_{\cU}(\psi_{\beta})$. Replacing $\psi_{\beta}$ by its convex envelope $\conve_{\cU}(\psi_{\beta})$, we obtain the convex outer approximation of $ \cSt$ in \eqref{eq.st2}:
\begin{equation}
\label{relcss}
  \bcSt\deq\{(u,v)\in \cU \times \cV: \conve_{\cU}(\psi_{\beta})(u) \le \psi_{\gamma}(v)\}.
\end{equation}

By using this extended formulation, our convex nonlinear relaxation of SP contains additional auxiliary variables. In particular, we need $2^h$ variables $\lambda_q$ to represent each convex envelope. For most SP problems in \minlplib where the degrees of the signomial terms are less than 6 and $h$ is less than 3, the convex nonlinear relaxation is computationally tractable.

\subsection{Outer approximation cuts}

The extended formulation is not useful, so we propose a cutting plane algorithm to separate valid linear inequalities in $(u,v)$-space from the extended formulation of the convex outer approximation. This algorithm generates a low-dimensional projected approximation of $\bcSt$. Moreover, the projection procedure converts the convex nonlinear relaxation into an LP relaxation, which is suitable for many solvers.

Given a point $(\relx{u}, \relx{v}) \in \cU \times \cV$, the algorithm determines whether it belongs to $\bcSt$. This verification can be done by checking the sign of $\conve_{\cU}(\psi_{\beta})(\relx{u}) - \psi_{\gamma}(\relx{v})$.  If $\conve_{\cU}(\psi_{\beta})(\relx{u}) - \psi_{\gamma}(\relx{v}) \le 0$, then $(\relx{u}, \relx{v}) \in \bcSt$. 

Since $\conve_{\cU}(\psi_{\beta})$ is a convex polyhedral function, our cutting plane algorithm evaluates the function by searching for an affine underestimator $a \cdot u+b$ of $\conve_{\cU}(u)$ such that $a \cdot \relx{u}+b = \conve_{\cU}(\relx{u})$, which is achieved by underestimating algorithms. If $ (\relx{u}, \relx{v}) \notin \bcSt$, then   $a \cdot u+b \le  \psi_{\gamma}(v)$ is   a valid nonlinear inequality of $\bcSt$. Subsequently, our cutting plane algorithm linearizes this inequality, resulting in  an outer approximation cut $a \cdot u+b \le  \lin{\psi_{\gamma}}{\relx{v}}(v)$: we recall that $\lin{\psi_{\gamma}}{\relx{v}}(v)$ is the linearization of $ \psi_{\gamma}(v)$ at $\relx{v}$  defined in Eq.~\eqref{def.Xi}.

We present our first \emph{LP-based} underestimating algorithm, which is used in our experiments.
 Due to \Cref{thm.poly}, we can solve the following LP to find the affine underestimator:
\begin{equation}
    \max_{a \in \bR^{h}, b \in \bR} a \cdot \relx{u}+b \quad \suc \forall q \in Q \; a \cdot q+b \le \psi_{\gamma}(q),
\end{equation}
where we omit the linear constraints that  bound $(a,b)$. The maximum value resulting from this LP is exactly $\conve_{\cU}(\psi_{\beta})(\relx{u})$. The affine underestimator $a \cdot u+b$ is called an \emph{facet} of the envelope $\conve_{\cU}(\psi_{\beta})$, if $a \cdot u +b \le t $ is a facet of $\epi(\conve_{\cU}(\psi_{\beta}))$. We note that the solution of the LP is not necessarily a facet, and the number of constraints is $2^h$.

We next give another enumeration-based underestimating algorithm. As $\psi_\beta$ is also concave, we recall the  characterization \cite{tawarmalani13} of the convex envelopes of concave functions $f$ over hypercubes.  A set of $h$-dimensional polyhedra $P_1 , \dots, P_t \subseteq \cU$ forms a \emph{triangulation} (\ie \emph{simplicial covers}) of $\cU$, if: (i) $\cU = \cup_{i \in [t]} P_t$; (ii) $P_i \cap P_j$ is a (possibly empty) face of both $P_i$ and $P_j$; (iii) each $P_i$ is an ($h$-)simplex. This means that each $P_i$ is  the convex hull of $h+1$ affine independent points (denoted as $S_i$).  We restrict our interests in triangulations that \emph{do not add vertices}, \ie every $S_i$ is a subset of the  vertices $Q$ of $\cU$.  We know that  an appropriate triangulation gives the convex envelope of $f$.

\begin{lemma}\cite[Thm.~2.4]{tawarmalani13} \label{lem.tri1}
For any concave function $f$, there exists a triangulation $\{P_i\}_{i \in [t]}$ of $\cU$ such that the convex envelope of $f$ over $\cU$ can be computed by interpolating $f$ affinely over each simplex $P_i$.
\end{lemma}

 However, it is non-trivial to find such an ``appropriate'' triangulation.
  To explain \Cref{lem.tri1}, any set $S   \deq  \{u^1,\dots,u^{h+1}\} \subseteq Q$ of $h+1$ affine independent points determines a function over $\bR^h$ via the following affine combination:
\begin{equation}
\label{eq.interpfunc}
	f_S(u) \deq \left \{\sum_{j \in [h+1]}\lambda_j f(u^j) : \exists \lambda \in \bR^{h+1} \, \sum_{j \in [h+1]}\lambda_j=1 \land \sum_{j \in [h+1]}\lambda_j u^j = u\right \}.
\end{equation}
Because of the affine independence of $S$, the barycentric coordinate $\lambda$ is unique for any $w$ in the above affine combination. We can consider $f_S$ as a single-valued affine function and call it the \emph{interpolation function} induced by $S$. Since $f_S$ interpolates $f$ at $S$, we can solve the linear system $a \cdot u + b = f(u)$ (for  $u \in S$) to compute $a,b$ that define $f_S$.  It follows from that \cite[Cor. 2.6]{tawarmalani13}, if $f_S$ underestimates $f$ at any point of $Q$, then $f_S$ is a facet of $\conve_{\cU}(f)$. We call such an $S$ \emph{facet-inducing}.

This result implies that we can focus on $h$-simplices instead of triangulations, since we want to find an affine underestimator for $f = \psi_\beta$. Our enumeration-based underestimating algorithm  finds the set of $h+1$ affine independent points in $Q$ such that the interpolation function $f_S$ is an underestimator of $f$. The algorithm outputs the greatest interpolation function at the point $\relx{u}$.

Finally, we explore another property of $ \psi_\beta$ that may help us reduce the search space.
 To simplify our representation, we translate and scale the domain of $\psi_{\beta}$ to $[0,1]^h$. This leads to a new function $s(w) \deq \psi_{\beta}(u)$, where for all $j \in [h]$, $u_j \deq \overline{u}_j + (\overline{u}_j - \underline{u}_j)w_j$. After these transformations, $\relx{u}$ becomes $\relx{w}$, the transformed domain $\cU$  of $u$ becomes $[0,1]^h$, and we denote the set of its vertices by the binary hypercube $Q = \{0,1\}^h$. W.l.o.g., we focus on the study and computation of facets of $\conve_{\cU}(s)$.

A set $D \subseteq \bR^h$ is called a \emph{product set}, if  $D = \bigtimes_{j \in [h]}D_j$ for $D_j \subseteq \bR$. Moreover, a function $f: D  \to \bR$ is \emph{supermodular} over $D$ (\cite[Sec. 2.6.1]{topkis2011supermodularity}), if the increasing difference condition holds:
for all $w^1,w^2 \in D, d \in \bR^h_{+}$ such that $w^1 \le w^2 $ and $ w^1+d, w^2+d \in D$, $f(w^1+d) -f(w^1)\le f(w^2 + d) -f(w^2)$. We find that  
the following operations preserve supermodularity. 

\begin{lemma}
\label{lem.subtrans}
    Let $w' \in \bR^h, \rho \in \bR^h_{++}$,  and let $D'$ be a product subset of $D$. The following results hold: (restriction) $f$ is supermodular over $D'$;(translation) $ f(w+w')$ is supermodular over $D-d$; (scaling) $f(\rho * w)$ is supermodular over $D/\rho$, where $+, -, *, /$ are taken entry-wise.
\end{lemma}
\begin{proof}
    The results follow from the definition.
\end{proof}

 We note that when $D = Q = \{0,1\}^h$, $d$ is in $Q$.  We observe a useful property of $s$. 
 
\begin{proposition}
\label{lem.sup}
The function $s$ is supermodular over $Q$ (\ie $\{0,1\}^h$). Moreover, $\conve_{\cU}(s) = \conve_{\cQ}(s) $.
\end{proposition}
\begin{proof} 
According to \cite[Example 2.6.2]{topkis2011supermodularity}, the signomial term
$\psi_\alpha$ with $\alpha > 0$  is supermodular  over $\bR^h_{+}$. This implies that the power function $\psi_\beta$ is supermodular over $\bR^h_+$. By \Cref{lem.subtrans}, $s$ is supermodular over $\cU = [0,1]^h$. As $Q=\{0,1\}^h$ is a product subset of $\cU$,   $s$  is supermodular over $Q$. After the scaling and translation, $s$ is still concave. By \Cref{thm.poly}, $\conve_{\cU}(s) = \conve_{\cQ}(s)$.
\end{proof}

Finding facets of $s$ could be reduced to a more genera problem of finding facets of supermodular functions over binary hypercubes. We note that a similar argument can show that both power functions and multilinear terms over any product subset of $\bR^h_+$ are supermodular.  

One may exploit  the increasing difference property to determine candidate sets of affine independent points when searching for facets. When
 $h = 2,$ we provide explicit projected formulations of convex envelopes of power functions. As a result, our cutting plane algorithm can efficiently separate outer approximation cuts for low-order problems.  For $h = 1$, the only facet is $s(0) + (s(1)- s(0))w_1$. 
\subsection{Projected convex envelopes in the bivariate case}

We present a general  characterization of projected  convex envelopes of supermodular functions  $f$ that is a restriction of a concave function. This gives a closed-form expression of the convex envelope of $s$ in the bivariate case. We can use a bit representation to denote binary points in $\{0,1\}^2$. For example, $10$ denotes the point $w$ that $w_1 = 1$ and $w_2 = 0$. For an affine function $a\cdot w + b$, we call  binary points in $\{0,1\}^2$ where $a\cdot w + b$ equals $f(w)$  its \emph{interpolating points}.

 Using the above result, we can construct an envelope-inducing family for bivariate supermodular functions. Let 
\begin{equation}
\label{eq.induce2}
    S^2_1 \deq \{00,10,01\}, S^2_2 \deq 
     \{11,10,01\}.
\end{equation}

One can find that $\conv(S^2_1)=\{(w_1,w_2) \in [0,1]^2:w_1 + w_2 \le 1 \}, \conv(S^2_2)=\{(w_1,w_2) \in [0,1]^2: w_1 + w_2 \ge 1 \}$ are two triangles in $[0,1]^2$.  We have that 
     	\begin{align*}
    	& f_{S^2_1}(w) =  f(00) + (f(10) -f(00)) w_1 + (f(01) -f(00)) w_2,\\
    	& f_{S^2_2}(w) = f(11) + (f(01) -f(11)) (1 - w_1) + (f(10) -f(11)) (1 - w_2).
	\end{align*}
 We show that these two affine functions define the convex envelope of $f$.
 
 \begin{theorem}
   Given $f:[0,1]^2 \to \bR$ a concave function that has a supermodular restriction over $\{0,1\}^2$, $\{S^2_k\}_{k \in [2]}$ as in \eqref{eq.induce2}  gives a triangulation of  $[0,1]^2$ and induce facets of $\conve_{[0,1]^2}(f)$.
 \end{theorem}
 \begin{proof}
 We know that $\conve_{[0,1]^2}(f) = \conve_{\{0,1\}^2}(f) $.
 It is easy to see that, for all $k \in [2]$, $S^2_k$ is affinely independent and $\{\conv(S^2_k)\}_{k \in [2]}$ is a triangulation of $[0,1]^2$.  Therefore, it suffices to show that $\{S^2_k\}_{k \in [2]}$  is facet-inducing, \ie $f_{S^2_1}, f_{S^2_2}$ are affine underestimators of $f$. 
 
 \textbf{Case i.}
 We note that, for all $w \in S^2_1 = \{00,10,01\}$, $f_{S^2_1}(w) = f(w).$ Note that $\{0,1\}^2\smallsetminus S^2_1 = \{11\}.$ It follows from the definition of the affine function $f_{S^2_1}$  that
 \begin{equation*}
     f_{S^2_1}(11)= f_{S^2_1}(10) + (f_{S^2_1}(01) - f_{S^2_1}(00)) = f(10) + (f(01) - f(00)).
 \end{equation*}
 It follows from the supermodularity of $f$ that 
 \begin{equation*}
     f(10) + (f(01) - f(00)) \le  f(10) + (f(11) - f(10)) = f(11).
 \end{equation*}
 Thereby, $f_{S^2_1}$ underestimates $f$.

  \textbf{Case ii.}
 We note that, for all $w \in S^2_2 = \{11,10,01\}$, $f_{S^2_2}(w) = f(w).$  Note that $\{0,1\}^2\smallsetminus S^2_2 = \{00\}.$ It follows from the definition of the affine function $f_{S^2_2}$  that
 \begin{equation*}
     f_{S^2_2}(00)= f_{S^2_2}(10) - (f_{S^2_1}(11) - f_{S^2_1}(01)) = f(10) + (f(11) - f(01)).
 \end{equation*}
 It follows from the supermodularity of $f$ that 
 \begin{equation*}
     f(10) - (f(11) - f(01)) \le  f(10) - (f(10) - f(00)) = f(00),
 \end{equation*}
which concludes the proof.
 \end{proof}

\subsection{Alternative convex outer approximations}
\label{sec.alter}
According to \Cref{sec.siliftfree}, we can have infinitely many DCC formulations of $\cSt$
parametrized by a scalar $\theta$:
\begin{equation*}
	\cSt^\theta \deq  \{(u, v) \in \cU \times \cV:\, \psi_{\theta \beta}(u) - \psi_{ \theta \gamma}(v) \le 0\},
\end{equation*}
where  $ \max(\lVert \beta \rVert_1, \lVert \gamma \rVert_1 ) = 1$ and $0 < \theta  \le 1$. Notice that $\cSt^1$ is used to construct the convex outer approximation  of $\cSt$. Alternatively, we have other convex outer approximations derived from $\cSt^\theta$:
\begin{equation*}
	\bcStt{\theta} \deq  \{(u, v) \in \cU \times \cV:\,   \conve_{\cU}(\psi_{\theta\beta})(u) - \psi_{ \theta \gamma}(v) \le 0\}.
\end{equation*}
For any $\theta, \theta' \in (0,1]$, $\cSt^\theta = \cSt^{\theta'}$, but $\bcStt{\theta} $ could be different from $\bcStt{\theta'}$. To generate the tightest outer approximation cuts, one may ask which $\theta$ yields the smallest convex outer approximation $\bcStt{\theta}$. We show that $\theta = 1$ is optimal in this sense.

We express $\cSt^\theta$ as follows:
\begin{equation*}
	\bcStt{\theta} \deq  \{(u, v) \in \cU \times \cV:\,   (\conve_{\cU}(\psi_{\theta\beta})(u))^{1/\theta} \le \psi_{ \gamma}(v) \}.
\end{equation*}
Since the right hand side $\psi_{ \gamma}(v)$ of the inequality does not depend on $\theta$, we check the value of the left hand side $(\conve_{\cU}(\psi_{\theta\beta})(u))^{1/\theta}$ at every point $u \in \cU$. We have the following observation on the bound of $(\conve_{\cU}(\psi_{\theta\beta})(u))^{1/\theta}$.

\begin{proposition}
    Given $u \in \cU$, for any $\theta \in (0,1]$, $(\conve_{\cU}(\psi_{\theta\beta})(u))^{1/\theta} $ is not greater than $\conve_{\cU}(\psi_{\beta})(u)$.
\end{proposition}
\begin{proof}
    According to \Cref{lem.tri1}, $\conve_{\cU}(\psi_{\beta})(u) = f_S(u)$, where $S$ is the set of $h+1$ affine independent points $u^j$ in the vertices $Q$ of $\cU$, and  the interpolation function $f(u)$ is taken as $\psi_{\beta}(u)$. Given the combination form \eqref{eq.interpfunc} of $f_S$, we express $\conve_{\cU}(\psi_{\beta})(u) = f_S(u) = \sum_{j \in [h+1]} \lambda_j \psi_{\beta}(u^j)$. Note that all $\lambda_j \ge 0$ (because $u \in \cU$), thus, the expression is indeed a convex combination form. Due to \Cref{thm.poly}, $\conve_{\cU}(\psi_{\theta\beta})(u)$ is the minimum of all convex combinations  $\sum_{q \in Q}\lambda_q\psi_{\theta\beta}(q)$. Thus,  $\conve_{\cU}(\psi_{\theta\beta})(u)$ is at most the particular convex combination  $\sum_{j \in [h+1]}\lambda_j\psi_{\theta\beta}(u_j)$. As $1/ \theta \ge 1$, $t^{1/\theta}$ is convex and non-decreasing w.r.t. the  indeterminate $t$. It follows from  the Jensen's inequality of convex function that \[(\conve_{\cU}(\psi_{\theta\beta})(u))^{1/\theta}  \le (\sum_{j \in [h+1]} \lambda_j \psi_{\beta}(u^j)^\theta)^{1/\theta}\le \sum_{j \in [h+1]} \lambda_j \psi_{\beta}(u^j),\]
    where the last convex combination is exactly $ \conve_{\cU}(\psi_{\beta})(u)$.
 \end{proof}

We then arrive at the conclusion on the optimality of $\theta = 1$.
\begin{corollary}
   For any $\theta \in (0,1]$, $\bcSt = \bcStt{1} \subseteq \bcStt{\theta} $.
\end{corollary}
\begin{proof}
    It is because $(\conve_{\cU}(\psi_{\theta\beta})(u))^{1/\theta}$ underestimates $ \conve_{\cU}(\psi_{\beta})(u)$.
\end{proof}

This explains why we choose $\theta = 1$ for our DCC formulation.
Note that the convex outer approximation derived from this formulation   may not be the convex hull of  $\cSt$.

 \subsection{Convexity and reverse-convexity}
 \label{sec.conv}
Our cutting plane algorithm can detect convexity/reverse-convexity of signomial-term sets. The detection is easily done by  the normalized DCC formulation, which gives another advantage.

Denote by $e^k_j$ and $e^h_j$ the $j$-th unit vector in $\bR^h$ and $\bR^k$, respectively. Then, we have the following observations:
\begin{enumerate}
	\item[i)] if $\norm{\beta}_1 = 1,\gamma = 0$, \ie $\psi_{\beta}$ is concave and $\psi_{\gamma}$ is 1, then $\cSt$ is reverse-convex;
	\item[ii)] if $\norm{\beta}_1 \le 1,\gamma = e^k_j$ for some $j \in [k]$, \ie $\psi_{\beta}$ is concave and $\psi_{\gamma}$ is a linear univariate function, then $\cSt$ is reverse-convex;
   \item[iii)] if $\beta = e^h_j, \norm{\gamma}_1 \le 1$ for some $j \in [h]$, \ie
$\psi_{\beta}$ is a linear univariate function and $\psi_{\gamma}$ is concave, then  $\cSt$ is convex;
   \item[iv)] if $\norm{\beta}_1 = 0, \norm{\gamma}_1 = 1$, \ie $\psi_{\beta}$ is 1  and $\psi_{\gamma}$ is concave, then $\cSt$ is convex.
\end{enumerate}

We note that similar results are found in \cite{chen2009note,maranas1995finding}. The results in  \cite{chen2009note} are proved by checking the negative/positive-semidefiniteness of the Hessian matrix of a signomial term. According to the normalized DCC formulation, the results are evident.

\section{Computational results}
\label{sec.comp}

In this section, we conduct computational experiments to assess the efficiency of the proposed valid inequalities.

The \minlplib dataset includes instances of MINLP problems containing signomial terms, and some of these instances are SP problems. To construct our benchmark, we select instances from \minlplib that satisfy the following criteria: (i) the instance contains signomial functions or polynomial functions, (ii) the continuous relaxation of the instance is nonconvex. Our benchmark consists of a diverse set of 251 instances in which nonlinear functions consist of signomial and other functions. These problems come from practical applications and can be solved by general purpose solvers.

Experiments are performed on a server with Intel Xeon W-2245 CPU @ 3.90GHz, 126GB main memory and Ubuntu 18.04 system. We use \scip 8.0.3 \cite{bestuzheva2023global} as a framework for reading and solving problems as well as performing cut separation. \scip is integrated with \texttt{CPLEX} 22.1 as LP solver and \texttt{IPOPT} 3.14.7 as NLP solver.

We evaluate the efficiency of the proposed valid inequalities in four different settings. In the first setting, denoted \disable, none of the proposed valid inequalities is applied. In the second setting, denoted \oc, only the outer approximation cuts are applied. The third setting, denoted \ic, applies only to the intersection cuts. The fourth setting combines both the \oc and \ic settings by applying both cuts. We let \scip's default internal cuts handle univariate signomial terms and multilinear terms. Our valid inequalities only handle the other high-order signomial terms. The source code, data, and detailed results can be found in our online repository: \href{https://github.com/lidingxu/ESPCuts}{github.com/lidingxu/ESPCuts}.

 Each test run uses \scip with a particular setting to resolve an instance. To solve the instances, we use the \scip solver with its sBB algorithm and set a time limit of 3600 seconds. In our benchmark, there are 150 instances classified as \emph{affected} in which at least one of the settings \oc, \ic, and \oic settings adds cuts. Among the affected instances, there are 86 instances where the default \scip configuration (\ie \disable setting) runs for at least 500 seconds. Such instances are classified as \emph{affected-hard}. For each test run, we measure the runtime, the number of sBB search nodes, and the relative open duality gap.

To aggregate the performance metrics for a given setting, we compute shifted geometric means (SGMs) over our test set. The SGM for runtime includes a shift of 1 second. The SGM for the number of nodes includes a shift of 100 nodes. The SGM for relative gap includes a shift of 1$\%$. We also compute the SGMs of the performance metrics over the subset of affected and affected-hard instances. The performance results are shown in \Cref{tb.perf}, where we also compute the relative values of the SGMs of the performance metrics compared to the \disable setting. Our following analysis is based on the results of the affected and affected-hard instances. Moreover, we display the absolute value of  the averaged separation time versus the absolute value of the averaged total runtime of each setting. We find that the separation time is much shorter than  the total runtime.

\begin{table} [htbp]
\centering
\resizebox{\columnwidth}{!}{
\addtolength{\tabcolsep}{-0.28em}
\begin{tabular}{cc|*{4}{c}|*{4}{c}|*{4}{c}}
\toprule
\multicolumn{2}{c|}{\multirow{2}{*}{Setting}}  &
\multicolumn{4}{c|}{All (\#251)}    &
\multicolumn{4}{c|}{Affected (\#150)} &
\multicolumn{4}{c}{Affected-hard (\#86)}  \\
& &
{solved} &
{nodes}     &
{time} &
{gap}     &
{solved} &
{nodes}     &
{time} &
{gap}   &
{solved} &
{nodes}     &
{time} &
{gap}   \\
\midrule
\multirow{2}{*}{\disable} & absolute  & \multirow{2}{*}{138} &  6510 & 0/122 & 4.7\%  & \multirow{2}{*}{71} &  15592 & 0/253 & 5.7\% & \multirow{2}{*}{7} & 175973 & 0/3600 & 26.7\%  \\ 
 & relative &  & 1.0 & 1.0 & 1.0 & & 1.0 & 1.0 & 1.0 & & 1.0 & 1.0 & 1.0 \\
 \midrule
\multirow{2}{*}{\oc} & absolute  & \multirow{2}{*}{140}  &  5954 & 1/118 & 4.5\%  & \multirow{2}{*}{73} & 13443 & 2/241 & 5.4\% &  \multirow{2}{*}{10} & 115262 & 9/2872 & 23.3\%   \\ 
 & relative &  &    0.91 & 0.97 & 0.97 & & 0.86 & 0.95 & 0.95 & & 0.65 & 0.8 & 0.87 \\
 \midrule
 \multirow{2}{*}{\ic} & absolute  & \multirow{2}{*}{140} & 6144 & 1/122 & 4.4\% & \multirow{2}{*}{73} &  14081 & 2/252 & 5.2\% & \multirow{2}{*}{10} & 128072 & 5/2994 & 22.0\%  \\ 
 & relative &  & 0.94 & 1.0 & 0.95 & &  0.9 & 0.99 & 0.91 &  & 0.73 & 0.83 & 0.82\\
 \midrule
 \multirow{2}{*}{\oic} & absolute  & \multirow{2}{*}{139} &  5934 & 1/117 & 4.6\% &  \multirow{2}{*}{72} & 13275 & 3/236 & 5.6\% & \multirow{2}{*}{10}  & 118054 & 10/2758 & 23.0\%\\ 
 & relative &  &  0.91 & 0.96 & 0.99 & & 0.85 & 0.93 & 0.98 &  & 0.67 & 0.77 & 0.86\\
\bottomrule
\end{tabular}
}
\caption{Summary of performance metrics on \minlplib instances.}\label{tb.perf}
\end{table}

First, we note that the proposed valid inequalities lead to the successful solution of 2 additional instances compared to the \disable setting. The \oc setting solves 2 more instances than the \disable setting.

The reductions in runtime and relative gap achieved by the \oc setting are 5\% and 5\%, respectively, for affected instances and 20\% and 13\%, respectively, for affected-hard instances. The \ic setting solves 2 more instances than the \disable setting. The reduction in runtime and relative gap achieved by the \ic setting is 1\% and 9\% for affected instances and 17\% and 14\% for affected-hard instances, respectively. The \oic setting resolves 1 additional instance compared to the \disable setting. The reduction in runtime and relative distance achieved by the \oic setting is 7\% and 2\%, respectively, for affected instances and 23\% and 14\%, respectively, for affected-hard instances.

We note that the runtime does not provide much information about affected-hard instances, since only 10 instances can be solved within 3600 seconds. For these instances, the gap reduction is more useful to measure the reduction of the search space by the proposed valid inequalities. However, for all affected instances, the runtime is still important because it measures the speedup due to the valid inequalities.

Second, we find that all cut settings have a positive effect on \scip performance, although the magnitude of the reduction varies. When we compare the \oc and \ic settings, we find that the \oc setting leads to a larger reduction in runtime. This difference in runtime is due to the fact that computing intersection cuts requires extracting a simplified cone from the LP relaxation and applying bisection search along each ray of the cone. These procedures require more computational resources compared to the construction of outer approximation cuts.

On the other hand, the \ic setting shows better performance in terms of reducing gaps. Intersection cuts approximate the intersection of a signomial-term set with the simplicial cone, while outer approximation cuts approximate the intersection of a signomial-term set with a hypercube. Around the relaxation point, the simplicial cone usually provides a better approximation than the hypercube. Therefore, \ic achieves a greater reduction in the relative gap. However, the better simplicial conic approximation does yield a significant improvement compared to the hypercubic approximation.

Finally, the \oic setting combines both the \oc and \ic settings and achieves the best reduction in runtime. However, for affected and affected-hard instances, the setting shows different gap reduction results. In fact, the results for affected-hard instances give more insight, since the goal of the valid inequalities is to speed up convergence for hard instances. In this sense, the \oic setting achieves almost the best result, so it carries the best of both valid inequalities. However, the improvement compared to each setting is not significant.

\begin{figure}[h]
	\centering
	\includegraphics[width=0.9\textwidth]{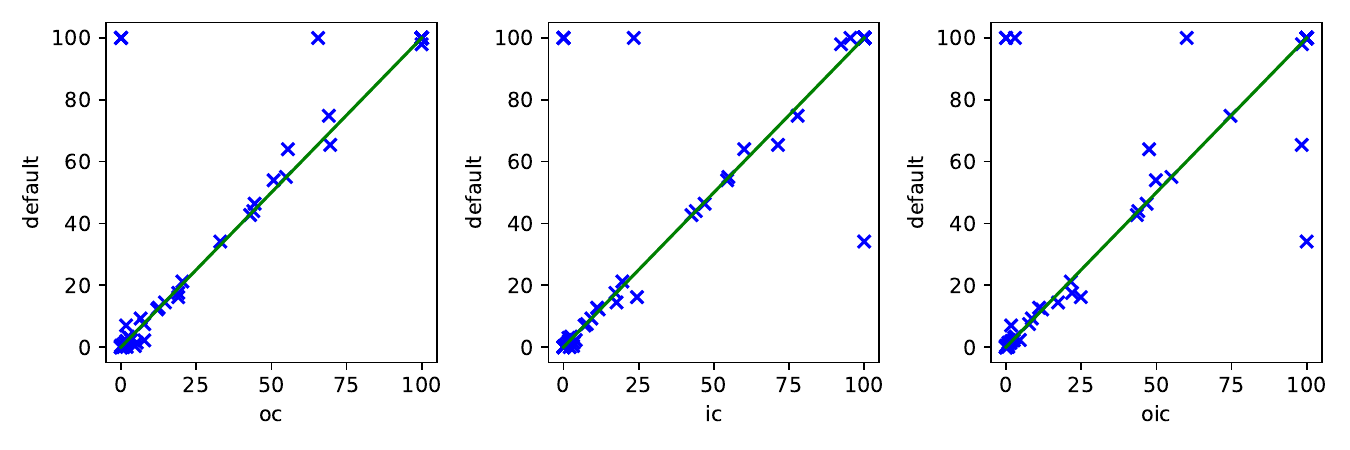}
	\caption{Relative gaps (in percentage) between pairs of settings for affected and unsolved (by the \disable) instances}
	\label{fig:comp2}
\end{figure}

We next look at instance-wise results on affected instances that are not solved by the \disable setting. The scatter plots in \Cref{fig:comp2} compare the  relative gaps of such instances obtained by different settings. We find that, many data  points (of gaps less than 20\%) are around the diagonal line, and these unbiased results mean that they are not affected much by cutting planes. However, there are some data points (of gaps more than 40\%) above the diagonal line,  especially noticing those far in the top, so cutting planes achieve much smaller gaps than the \disable setting on these hard  instances.

In summary, the performances of the \oc and \ic settings are comparable. They can lead to smaller duality gaps, which is desirable for solvers, and one can use either of them. Moreover, the combination of both cuts enhances performance slightly.

\section{Conclusion and discussions}

In this paper we study valid inequalities for SP problems and propose two types of valid linear inequalities: intersection cuts and outer approximation cuts. Both are derived from normalized DCC formulations of signomial-term sets. First, we study general conditions for maximal $\cS$-free sets. We construct maximal signomial term-free sets from which we generate   intersection cuts. Second, we construct convex outer approximations of signomial-term sets within hypercubes. We provide extended formulations for the convex envelopes of concave functions in the normalized DCC formulations. Then we separate valid inequalities for the convex outer approximations by projection. Moreover, when $h=2$, we use supermodularity to derive a closed-form expression for the convex envelopes.

We present a comparative analysis of the computational results obtained with the \minlplib instances. This analysis demonstrates the effectiveness of the proposed valid inequalities. The results show that intersection cuts and outer approximation cuts have similar performance and their combination takes the best of each setting. In particular, it is easy to implement outer approximation cuts in general purpose solvers.

 In the following, we have some further discussions that lead to some open questions and possible extensions of the proposed cutting plane algorithms.

\subsection{Signomial aggregation}

We currently deal with signomial terms explicitly present in the signomial terms, but our results can be extended to deal with multiple signomial terms. In the future, the proposed valid inequalities can approximate nonlinear aggregations of constraints that define the signomial lift. Specifically, given signomial constraints $\{\psi_{\alpha^i}(x) = y_i\}_{i \in [r]}$ with any exponent vector $\zeta \in \bR^r$, we can employ \emph{signomial aggregation} to generate a new signomial constraint: $\psi_{(\sum_{i \in [r]}\zeta_i\alpha^i)}(x) = \psi_{\zeta}(y)$. This constraint is valid for the signomial lift and encodes more variables and terms. Next, we can apply the DCC reformulation to the constraints $\psi_{(\sum_{i \in [r]}\zeta_i\alpha^i)}(x) \le \psi_{\zeta}(y)$ and $\psi_{(\sum_{i \in [r]}\zeta_i\alpha^i)}(x) \ge \psi_{\zeta}(y)$. Finally, we can separate the proposed valid inequalities. As far as we know, the signomial aggregation operator is not yet used for polynomial programming, since it outputs a signomial constraint.

\subsection{Signomial constraints}
Through lifting signomial terms,
we have studied the extended formulation of SP. The proposed methods could be used for relaxing  signomial constraints in the projected formulation of SP, but this may require a global transformation of variables. We can always write a signomial constraint as follows:
\begin{equation}
\label{eq.sigconsside}
    \sum_{i \in I_1} b_i \psi_{\alpha^i}(x) \le  \sum_{i \in I_2} b_i \psi_{\alpha^i}(x),
\end{equation}
where, for all $i \in I \deq I_1 \cup I_2$,  $b_i \ge 0$ and $\alpha_i \in \bR^n$. We want the signomial terms to have only positive exponents. As the both sides of the signomial constraint \eqref{eq.sigconsside} are non-negative, we can multiply both sides by a signomial term $\psi_{\alpha^0}(x)$ with $\alpha^0 \ge 0$, which should yield all  $\beta^i \deq \alpha^i + \alpha^0 \ge 0$. This  reformulates  the signomial constraint \eqref{eq.sigconsside} as follows:
\begin{equation}
\label{eq.sigconsside1}
    \sum_{i \in I_1} b_i \psi_{\beta^i}(x) \le  \sum_{i \in I_2} b_i \psi_{\beta^i}(x).
\end{equation}

Note that $ \psi_{\beta^i}(x)$ only have positive exponents, but they are not necessarily power functions. For the constraint in reformulated signomial-term set in \eqref{eq.form2}, we applies powers  on two signomial terms to rescale their exponents, and we obtain a DCC constraint. However, this power rescaling  generally does not produce a DCC reformulation of \eqref{eq.sigconsside1}, because the rescaled term $(  \sum_{i \in I_1} b_i \psi_{\beta^i}(x) )^{\mu}$ for $\mu > 0$ could be nonconvex. Instead, we can use power transformation to  overcome this difficulty. Given $ \gamma \in \bR^n_{++}$,  denote $z = (x_j^{\gamma_j})_j$, and we note that $\psi_{\beta^i}(x)  =   \psi_{\beta^i / \gamma}(z)$, where $/$ is taken entry-wise.   When all $\norm{\beta^i / \gamma}_1 \le 1$,  every $\psi_{\beta^i}(z)$ is a power function. Therefore, the signomial constraint \eqref{eq.sigconsside1} is equivalent to the following DCC constraint:
\begin{equation}
\label{eq.ssss}
    \sum_{i \in I_1} b_i \psi_{\beta^i/ \gamma}(z) \le  \sum_{i \in I_2} b_i \psi_{\beta^i / \gamma}(z).
\end{equation}

Note that the SP can have other signomial constraints in $x$, and this global power transformation reformulates SP  in the variable space of $z$. We should choose an appropriate parameter $\gamma$ that transforms all signomial constraints into DCC constraints in $z$ as well, and such a $\gamma$ should satisfy that $\norm{\beta / \gamma}_1 \le 1$ for all exponents $\beta$ appearing in the reformulated signomial constraints as \eqref{eq.sigconsside1}. Then, we could apply the proposed cutting planes on this space. However, it is not easy to implement this  global power transformation  in current solvers, or such a transformation does not exist for problems mixed with signomial terms and other nonlinear functions.  We pose some open problems here. Which $\gamma$  yields DCC constraints  that result in maximal $\cS$-free sets ($\cS$ is taken as the feasible set defined by the constraint \eqref{eq.ssss})? As 
\Cref{thm.max} requires at least one part of the DCC function to be positive homogeneous of degree 1, could we reuse \Cref{thm.max} to find $\gamma$? We conjecture that such a $\gamma$ does not exist in general, because we have to ensure several $\norm{\beta^i / \gamma} = 1$.

\section*{Acknowledgments}

We thank Felipe Serrano for sharing his Ph.D. thesis \cite{serranothesis} with us. We appreciate the discussion with Antonia Chmiela,
Mathieu Besançon, and Stefan Vigerske on \texttt{SCIP}. Special thanks go to Ksenia Bestuzheva for helping us integrate the outer approximation cuts in the version 9.0 \cite{BolusaniEtal2024} of \texttt{SCIP}.
We also thank two anonymous reviewers for many helpful comments, which significantly improve the readability of the article.  This publication was supported by the ``Integrated Urban Mobility'' chair, which is supported by L’X - \'Ecole Polytechnique and La Fondation de l’\'Ecole Polytechnique. Under no circumstances do the partners of the Chair assume any liability for the content of this publication, for which the author alone is responsible.

\section*{Declarations}
No conflicts of interest with the journal or funders.


\section*{Appendix}
\begin{proof}[Proof of \Cref{thm.decomp}]
It suffices to consider the case that $\cC$ is a maximal $\epi(g'_i)$-free set in $\cG^i$. W.l.o.g., we can assume that $\cC,  \cG^i$ are full-dimensional in $\bR^{I_i}$. Since $\epi(g'_i)$ includes $\gr(g'_i)$, $\cC$, as an $\epi(g'_i)$-free set,  is also $\gr(g'_i)$-free. First, we prove that  $\cC$ is a maximal $\gr(g'_i)$-free set in $\cG^i$. Assume, to aim at a contradiction, that $\cC'$ is a $\gr(g'_i)$-free set that $\cC \cap \cG^i \subsetneq \cC' \cap \cG^i$.  Suppose that $\epi(g'_i)  \cap \inter(\cC' \cap \cG^i)$ is not empty and contains $(x'_{J_i}, y'_i)$.     As $\cC$ is  $\epi(g'_i)$-free, there exists a point $(x_{J_i}, y_i) \in \inter(\cC \cap \cG^i) \subseteq \inter(\cC' \cap \cG^i)$ such that $(x_{J_i}, y_i) \in \hyp(g'_i)$.  It follows from the continuity of $g'_i$  that there exists a point $(x^\ast_{J_i}, y^\ast_i) \in \gr(g'_i)$ in the line segment joining $(x_{J_i}, y_i)$ and $(x'_{J_i}, y'_i) $. As $\inter(\cC' \cap \cG^i)$ is convex,  we have that $(x^\ast_{J_i}, y^\ast_i) \in \inter(\cC' \cap \cG^i)$, which leads to a contradiction to $\gr(g'_i)$-freeness of $\cC'$. Therefore, $\epi(g'_i)  \cap \inter(\cC' \cap \cG^i)$ must be empty, so $\cC' \cap \cG^i  \subseteq \hyp(g'_i)$. This means that $\cC'$ is also $\epi(g'_i)$-free. However, note that  $\cC \cap \cG^i \subsetneq \cC' \cap \cG^i$, this contradicts with the fact that $\cC$ is a maximal $\epi(g'_i)$-free set in $\cG^i$. Therefore, $\cC$  is a maximal $\gr(g'_i)$-free set in $\cG^i$. Secondly, we  prove that $\bar{\cC}$ is a maximal $\cSl$-free set in $\cG$.
Assume, to aim at a contradiction, that there exists an $\cSl$-free set $\bar{\cD}$ in $\cG$ such that $ \bar{\cC} \cap \cG \subsetneq \bar{\cD} \cap \cG$. We look at their orthogonal projections on  $\bR^{I_i}$. It follows from the decomposability  that $\cC \cap \cG^i = \cC \cap \proj_{\bR^{I_i}} (\cG) = \proj_{\bR^{I_i}}(\bar{\cC} \cap \cG) \subseteq \proj_{\bR^{I_i}}(\bar{\cD} \cap \cG)$. Denote $\cD \deq \cl(\proj_{\bR^{I_i}}(\bar{\cD} \cap \cG))$, which is a closed convex set in $\cG^i$. Since $\bar{\cC}= \cC \times \bR^{I^c_i}$, $\cD$  must strictly include $\cC \cap \cG^i$. Note that $\cD$ is $\gr(g'_i)$-free. Since $\cC$ is a maximal $\gr(g'_i)$-free set in $\cG^i$, this implies that $\cC \cap \cG^i = \cD$, which leads to a contradiction.
\end{proof}

	\begin{proof}[Proof of \Cref{prop.max}]
 Let $\cC$ be a set satisfying the hypothesis.  Suppose, to aim at a contradiction, that there exists a closed convex set $\cC^\ast$ such that $\cC \subsetneq \cC^\ast$  and $\cC^\ast$ is an inner approximation of $\cF$. Then, there must exist  an open ball $B$ such that  \( B \subseteq \cF \smallsetminus\cC\) and $B \subseteq  \cC^\ast$. Let $z^\ast$ be the center of $B$, so   $z^\ast \in  \inter(\cF \smallsetminus \cC)$. W.l.o.g., we let $\cC^\ast = \conv(\cC \cup \{z^\ast\})$, which is a closed convex inner approximation of $\cF$. Since $z^\ast \notin \cC$, by the hyperplane separation theorem, there exists  \({\lambda} \in \dom(\sigma_{\cC})\) such that
	\begin{equation}
    	\label{eq.r0}
	{\lambda}\cdot z^\ast > \sigma_{\cC}({\lambda}).
	\end{equation}
	For any such $\lambda$, by the hypothesis, there exists a point \(z' \in \bd(\cF) \cap \bd(\cC)\) such that
	\begin{equation}
	\label{eq.r1}
  	{\lambda} \cdot z'  = \sigma_{\cC}({\lambda}),  
	\end{equation}
	and $z'$ is an exposing point of $\cC$. We want to show that,  for any \({\lambda'} \in \dom(\sigma_{\cC^\ast})\), \({\lambda'}\cdot z' < 	\sigma_{\cC^\ast}({\lambda})\). We consider the following three cases. First, we consider the case $\lambda'=\lambda$.  Because \(z^\ast \in \cC^\ast\), by the definition of support functions, we have that
 	\begin{equation}
     	\label{eq.r2}
     	{\lambda} \cdot z^\ast \le \sup_{z \in \cC^\ast}{\lambda}\cdot z = \sigma_{\cC^\ast}({\lambda}).
 	\end{equation}
	It follows from \eqref{eq.r0}, \eqref{eq.r1}, and \eqref{eq.r2} that
	\begin{align}
	\label{eq.r4}
    	{\lambda} \cdot z'  =  \sigma_{\cC}({\lambda}) < {\lambda} \cdot z^\ast \le \sigma_{\cC^\ast}({\lambda}) = \sigma_{\cC^\ast}({\lambda'}).
	\end{align} Second, we consider the case $\lambda ' = \rho \lambda$ for some $\rho > 0$. Since $\sigma_{\cC^\ast}$ is positively homogeneous of degree 1, it follows from \eqref{eq.r4} that $\lambda' \cdot z' =  {\rho \lambda} \cdot z' < \rho \sigma_{\cC^\ast}({ \lambda}) = \sigma_{\cC^\ast}({\lambda'}).$ Last, we consider the case ${\lambda'} \in \dom(\sigma_{\cC^\ast}) \smallsetminus\{{\rho \lambda}\}_{\rho > 0}$. By \Cref{lem.supconv}, $\sigma_{\cC} \le \sigma_{\cC^\ast}$. 
	By the hypothesis that $z'$ is an exposing point of $\cC$, provided that $\lambda' \ne \rho \lambda$,  we have that
	$
    	{\lambda'}\cdot z' < \sigma_{\cC}({\lambda'}) \le	\sigma_{\cC^\ast}({\lambda'}).
	$
	In summary, we have proved that for any \({\lambda'} \in \dom(\sigma_{\cC^\ast})\), \({\lambda'}\cdot z' < 	\sigma_{\cC^\ast}({\lambda'})\). So by \Cref{lem.supconv}, \(z' \in \inter(\cC^\ast)\). We find that $z' \in \bd(\cF) \cap \inter(\cC^\ast)$. This finding means a point near  $z'$ exists, which is in $\cC^\ast$, but not in $\cF$. Hence, $\cC^\ast$ is not an inner approximation of $\cF$, which leads to a contradiction.
	\end{proof}

\begin{proof}[Proof of \Cref{lem.homo}]
Given $z \in \dom(f)$, $f(\rho z)= \rho^d f(z) $ is a real number for any  $\rho \in \bR_{++}$, so $\inter(\dom(f))$ is a cone.
Suppose that $f$ is positive homogeneous of degree 1. For any $z \in \dom(f)$,
$
	\lin{f}{\breve{z}}(z) = f(\breve{z}) + \nabla  {f}(\breve{z}) \cdot (z - \breve{z}) =  \nabla  {f}(\breve{z}) \cdot z,
$
where  the second equation follows from Euler's homogeneous function theorem: $f(\breve{z}) = \nabla {f}(\breve{z}) \cdot \breve{z}$. For any \( z = \rho \breve{z}\) with \(\rho \in \bR_{++}\),
$
	\lin{f}{\breve{z}}(z)  = \nabla  {f}(\breve{z})\cdot \rho \breve{z} = \rho  \lin{f}{\breve{z}}(\breve{z}) =\rho f( \breve{z}) = f(\rho \breve{z}),
$
where the first and second equations follow from the previous result, the third follows from that $\lin{f}{\breve{z}}$ has the same value as $f$ at $\breve{z}$, and the last equation follows from the homogeneity.
\end{proof}

\bibliographystyle{siamplain}

\bibliography{reference}

\end{document}